\documentclass[10pt]{amsart}
\numberwithin{equation}{section}

\usepackage{amssymb}
\usepackage{amscd}

\newtheorem{thm}{Theorem}[section]
\newtheorem{lem}[thm]{Lemma}
\newtheorem{constrlem}[thm]{Construction-Lemma}

\theoremstyle{definition}
\newtheorem{dfn}[thm]{Definition}
\newtheorem{ntn}[thm]{Notation} 
\newtheorem{rmk}[thm]{Remark}
\newtheorem{exa}[thm]{Example}

\newtheorem{conj}[thm]{Conjecture}

\newenvironment{pf}
  {\noindent{{\em Proof: }}}
  {\qed \medskip}

\def\quod{\hskip 0.5em\relax }

\newcommand{\Pb}{\mathbf{P}}
\newcommand{\Ab}{\mathbf{A}}
\newcommand{\abs}[1]{\lvert#1\rvert}
\newcommand{\Z}{\mathbb{Z}}
\newcommand{\Q}{\mathbb{Q}}
\newcommand{\N}{\mathbb{N}}
\newcommand{\Ql}{\mathbb{Q}_{\ell}}
\newcommand{\C}{\mathbb{C}}
\newcommand{\Aut}{\mathrm{Aut}}
\newcommand{\s}{\mathbb{S}}
\newcommand{\V}{\mathbb{V}}
\newcommand{\Vla}{\mathbb{V}_{\lambda}}

\newcommand{\Mm}{\mathcal M}
\newcommand{\Hh}{\mathcal H}
\newcommand{\Mmb}{\overline{\mathcal M}}

\newcommand{\Ll}{\mathbf L}
\newcommand{\Eul}{\mathbf{e}}

\newcommand{\Eulc}{\mathbf{e}_c}

\newcommand{\UU}{\hat{U}}
\newcommand{\bb}{\hat{b}}
\newcommand{\nn}{\mathbf{n}}
\newcommand{\rr}{\mathbf{r}}

\begin{document}
\pagestyle{plain}

\title[Equivariant counts of points of $\mathcal{H}_{g,n}$]{Equivariant counts of points of the moduli spaces of pointed hyperelliptic curves} 
\author{Jonas Bergstr\"om}
\address{Korteweg-de Vries Instituut, Universiteit van Amsterdam, Plantage Muidergracht 24, 1018 TV Amsterdam, The Netherlands.}
\email{o.l.j.bergstrom@uva.nl}

\begin{abstract} 
We consider the moduli space $\Hh_{g,n}$ of $n$-pointed smooth hyperelliptic curves of genus $g$. In order to get cohomological information we wish to make $\s_n$-equivariant counts of the numbers of points defined over finite fields of this moduli space. 
We find recurrence relations in the genus that these numbers fulfill. Thus, if we can make $\s_n$-equivariant counts of $\Hh_{g,n}$ for low genus, then we can do this for every genus. Information about curves of genus $0$ and $1$ is then found to be sufficient to compute the answers for $\Hh_{g,n}$ for all $g$ and for $n \leq 7$. These results are applied to the moduli spaces of stable curves of genus $2$ with up to $7$ points, and this gives us the $\s_n$-equivariant Galois (resp. Hodge) structure of their $\ell$-adic (resp. Betti) cohomology. 
\end{abstract}

\maketitle

\section{Introduction}
By virtue of the Lefschetz trace formula, counting points defined over finite fields of a space gives a way of finding information on its cohomology. In this article we wish to count points of the moduli space $\Hh_{g,n}$ of $n$-pointed smooth hyperelliptic curves of genus $g$. On this space we have an action of the symmetric group $\s_n$ by permuting the marked points of the curves. To take this action into account we will make \emph{$\s_n$-equivariant} counts of the numbers of points of $\Hh_{g,n}$ defined over finite fields. 

For every $n$ we will find simple recurrence relations in the genus, for the equivariant number of points of $\Hh_{g,n}$ defined over a finite field. Thus, if we can count these numbers for low genus, we will know the answer for every genus. The hyperelliptic curves will need to be separated according to whether the characteristic is odd or even and the respective recurrence relations will in some cases be different.

When the number of marked points is at most $7$ we use the fact that the base cases of the recurrence relations only involve the genus $0$ case, which is easily computed, and previously known $\s_n$-equivariant counts of points of $\Mm_{1,n}$, to get equivariant counts for every genus. If we consider the odd and even cases separately, then all these counts are polynomials when considered as functions of the number of elements of the finite field. For up to five points these polynomials do not depend upon the characteristic. But for six-pointed hyperelliptic curves there is a dependence, which appears for the first time for genus $3$. 

By the Lefschetz trace formula, the $\s_n$-equivariant count of points of $\Hh_{g,n}$ is equivalent to the trace of Frobenius on the $\ell$-adic $\s_n$-equivariant Euler characteristic of $\Hh_{g,n}$. But this information can also be formulated as traces of Frobenius on the Euler characteristic of some natural local systems $\Vla$ on $\Hh_g$. By Theorem 3.2 in \cite{Jonas2} we can use this connection to determine the Euler characteristic, evaluated in the Grothendieck group of absolute Galois modules, of all $\Vla$ on $\Hh_g \otimes \overline{\Q}$ of weight at most~$7$. These result are in agreement with the results on the ordinary Euler characteristic and the conjectures on the motivic Euler characteristic of $\Vla$ on $\Hh_3$ by Bini-van der Geer in \cite{Bini-Geer}, the ordinary Euler characteristic of $\Vla$ on $\Hh_{2}$ by Getzler in \cite{G-euler}, and the $\s_2$-equivariant cohomology of $\Hh_{g,2}$ for all $g \geq 2$ by Tommasi in \cite{Orsolathesis}. 

The moduli stack $\Mmb_{g,n}$ of stable $n$-pointed curves of genus $g$ is smooth and proper, which implies purity of the cohomology. If the $\s_n$-equivariant count of points of this space, when considered as a function of the number of elements of the finite field, gives a polynomial, then using the purity we can determine the $\s_n$-equivariant Galois (resp. Hodge) structure of its individual $\ell$-adic (resp. Betti) cohomology groups (see Theorem 3.4 in \cite{Mbar4} which is based on a result of van den Bogaart-Edixhoven in \cite{EB}). All curves of genus~$2$ are hyperelliptic and hence we can apply this theorem to  $\Mmb_{2,n}$ for all $n \leq 7$. These results on genus $2$ curves are all in agreement with the ones of Faber-van der Geer in \cite{Faber-Geer1} and \cite{Faber-Geer2}. Moreover, for $n \leq 3$ they were previously known by the work of Getzler in \cite[Section 8]{G-2}. 

\section*{Acknowledgements}
The method I shall use to count points of the moduli space of pointed hyperelliptic curves follows a suggestion by Nicholas M. Katz. I thank Bradley Brock for letting me read an early version of the article \cite{Brock}, and Institut Mittag-Leffler for support during the preparation of this article. I would also like to thank my thesis advisor Carel Faber. 

\section*{Outline} 
Let us give an outline of the paper, where $\star_{\cdot}$ denotes the section.
\begin{itemize}
\item[$\star_2$] In this section we define $\s_n$-equivariant counts of points of $\Hh_{g,n}$ over a finite field $k$, and we formulate the counts in terms of numbers $a_{\lambda}|_g$, which are connected to the $H^1$'s of the hyperelliptic curves. 
\item[$\star_3$] The hyperelliptic curves of genus $g$, in odd characteristic, are realized as degree $2$ covers of $\Pb^1$ given by square-free polynomials of degree $2g+2$ or $2g+1$. The numbers $a_{\lambda}|_g$ are then expressed in terms of these polynomials in equation \eqref{eq-afp}. The expression for $a_{\lambda}|_g$ is decomposed into parts denoted $u_g$, which are indexed by pairs of tuples of numbers $(\nn;\rr)$. The special cases of genus $0$ and $1$ are discussed in Section~\ref{sec-g01}. 
\item[$\star_4$] A recurrence relation is found for the numbers $u_g$ (Theorem \ref{thm-rec1}). The first step is to use the fact that any polynomial can be written uniquely as a monic square times a square-free one. This results in an equation which gives $U_g$ in terms of $u_{h}$ for $h$ less than or equal to $g$, where $U_g$ denotes the expression corresponding to $u_g$, but in terms of all polynomials instead of only the square-free ones. The second step is to use that, if $g$ is large enough, $U_g$ can be computed using a simple interpolation argument. 
\item[$\star_5$] The recurrence relations for the $u_g$'s are put together to form a linear recurrence relation for $a_{\lambda}|_g$, whose characteristic polynomial is given in Theorem~\ref{thm-chareq}.
\item[$\star_6$] It is shown how to compute $u_0$ for any pair $(\nn;\rr)$.
\item[$\star_7$] Information on the cases of genus 0 and 1 is used to compute, for all $g$, $u_g$ for tuples $(\nn;\rr)$ of degree at most $5$, and $a_{\lambda}|_g$ of weight at most $7$. 
\item[$\star_8$] The hyperelliptic curves are realized, in even characteristic, as pairs $(h,f)$ of polynomials fulfilling three conditions. The numbers $u_g$ and $U_g$ are then defined to correspond to the case of odd characteristic. 
\item[$\star_9$] In even characteristic, a recurrence relation is found for the numbers $u_g$ (Theorem \ref{thm-rec1even}). Lemmas \ref{lem-disjoint} and \ref{lem-cover} show that one can do something in even characteristic corresponding to uniquely writing a polynomial as a monic square times a square-free one in odd characteristic. This results in a relation between $U_g$ and $u_h$ for $h$ less than or equal to $g$. Then, as in odd characteristic, a simple interpolation argument is used to compute $U_g$ for $g$ large enough. 
\item[$\star_{10}$] The same amount of information as in Section~\ref{sec-results} is obtained in the case of even characteristic. It is noted that $a_{\lambda}|_g$ is independent of the characteristic for weight at most $5$ (Theorem \ref{thm-ind}). This does not continue to hold for weight $6$ where there is dependency for genus at least $3$ (see Example \ref{exa-a16even}).
\item[$\star_{11}$] The counts of points of the previous sections are used to get cohomological information. This is, in particular, applied to $\Mmb_{2,n}$ for $n \leq 7$.
\item[$\star_{12}$] In the first appendix, a more geometric interpretation is given of the information contained in all the numbers $u_g$ of at most a certain degree (see Lemma~\ref{lem-bc-u}).
\item[$\star_{13}$] In the second appendix, we find that for sufficiently large $g$ we can compute the Euler characteristic, with $\mathrm{Gal}(\overline{\Q}/\Q)$-structure, of the part of the cohomology of sufficiently high weight, of some local systems $\Vla$ on $\Hh_g$. We will also see that these results are, in a sense, \emph{stable} in $g$.
\end{itemize}

\section{Equivariant counts} \label{sec-equiv}
Let $k$ be a finite field with $q$ elements and denote by $k_m$ a degree $m$ extension. Define $H_{g,n}$ to be the coarse moduli space of $\Hh_{g,n} \otimes \bar k$ and let $F$ be the geometric Frobenius morphism. 

The purpose of this article is to make $\s_n$-equivariant counts of the number of points defined over $k$ of $H_{g,n}$. With this we mean a count, for each element $\sigma \in \s_n$, of the number of fixed points of $F \sigma$ acting on $H_{g,n}$. Note that these numbers only depend upon the cycle type $c(\sigma)$ of the permutation $\sigma$.

Define $\mathcal{R}_{\sigma}$ to be the category of hyperelliptic curves of genus $g$ that are defined over $k$ together with marked points $(p_1,\ldots,p_n)$ defined over $\bar{k}$ such that $(F \sigma) (p_i)=p_i$ for all $i$. Points of $H_{g,n}$ are isomorphism classes of $n$-pointed hyperelliptic curves of genus $g$ defined over $\bar k$. For any pointed curve $X$ that is a representative of a point in $H_{g,n}^{F \sigma}$, the set of fixed points of $F\sigma$ acting on $H_{g,n}$, there is an isomorphism from $X$ to the pointed curve $(F\sigma) X$. Using this isomorphism we can descend to an element of $\mathcal{R}_{\sigma}$ (see \cite[Lem.~10.7.5]{Katz-Sarnak}). Therefore, the number of $\bar k$-isomorphism classes of the category $\mathcal{R}_{\sigma}$ is equal to $\abs{H_{g,n}^{F \sigma}}$.

Fix an element $Y=(C,p_1,\ldots,p_n)$ in $\mathcal{R}_{\sigma}$. We then have the following equality (see \cite{Geer} or \cite{Katz-Sarnak}):
$$\sum_{\substack{[X] \in \mathcal{R}_{\sigma}/\cong{k} \\ X\cong_{\bar{k}}Y}} \frac{1}{\abs{\Aut_{k}(X)}} = 1.$$
This enables us to go from $\bar k$-isomorphism classes to $k$-isomorphism classes:
$$\abs{H_{g,n}^{F \sigma}} = \sum_{[Y] \in \mathcal{R}_{\sigma}/\cong_{\bar{k}}} 1 =  \sum_{[Y] \in \mathcal{R}_{\sigma}/\cong_{\bar{k}}} \sum_{\substack{[X] \in \mathcal{R}_{\sigma}/\cong{k} \\ X\cong_{\bar{k}}Y}} \frac{1}{\abs{\Aut_{k}(X)}}  = \sum_{[X] \in \mathcal{R}_{\sigma}/\cong_{k}} \frac{1}{\abs{\Aut_{k}(X)}}.$$

For any curve $C$ over $k$, define $C\bigl(\sigma\bigr)$ to be the set of $n$-tuples of distinct points $(p_1,\ldots,p_n)$ in $C(\bar{k})$ that fulfill $(F\sigma) (p_i)=p_i$. 
\begin{ntn} A partition $\lambda$ of an integer $m$ consists of a sequence of non-negative integers $\lambda_1,\ldots, \lambda_{\nu}$ such that $\abs{\lambda}:=\sum_{i=1}^{\nu} i \lambda_i=m$. We will write $\lambda=[1^{\lambda_1},\ldots,\nu^{\lambda_{\nu}}]$.
\end{ntn}

Say that $\tau \in \s_n$ consists of one $n$-cycle. The elements of $C\bigl(\tau \bigr)$ are then given by the choice of $p_1 \in C(k_n)$ such that $p_1 \notin C(k_i)$ for every $i<n$. By an inclusion-exclusion argument it is then straightforward to show that
$$\abs{C\bigl(\tau \bigr)}=\sum_{d | n} \mu(n/d)\,\abs{C(k_d)},$$
where $\mu$ is the M\"obius function. Say that $\lambda$ is any partition and that $\sigma \in \s_{\abs{\lambda}}$ has the property $c(\sigma)=\lambda$. Since $C\bigl(\sigma\bigr)$ consists of tuples of distinct points it directly follows that
\begin{equation} \label{eq-sigma}
\abs{C\bigl(\sigma \bigr)}=\prod_{i=1}^{\nu} \prod_{j=0}^{\lambda_i-1} \Biggl(\sum_{d | i} \Bigl(\mu(i/d)\,\abs{C(k_d)} -j  i \Bigr)\Biggr).
\end{equation}

Fix a curve $C$ over $k$ and let $X_1, \ldots, X_m$ be representatives of the distinct $k$-isomorphism classes of the subcategory of $\mathcal{R}_{\sigma}$ of elements $(D,q_1,\ldots,q_n)$ where $D \cong_k C$. For each $X_i$ we can act with $\Aut_k(C)$ which gives an orbit lying in $\mathcal{R}_{\sigma}$ and where the stabilizer of $X_i$ is equal to $\Aut_k(X_i)$. Together the orbits of $X_1, \ldots, X_m$ will contain $\abs{C\bigl(\sigma\bigr)}$ elements and hence we obtain
\begin{equation} \label{eq-equiv}
\abs{H_{g,n}^{F  \sigma}} = \sum_{[X] \in \mathcal{R}_{\sigma}/\cong_{k}} \frac{1}{\abs{\Aut_{k}(X)}} = \sum_{[C] \in \Hh_g(k)/\cong_k} \frac{\abs{C\bigl(\sigma\bigr)}}{\abs{\Aut_k(C)}}.  
\end{equation}

We will compute slightly different numbers than $\abs{H_{g,n}^{F  \sigma}}$, but which contain equivalent information. Let $C$ be a curve defined over $k$. The Lefschetz trace formula tells us that for all $m \geq 1$, 
\begin{equation} \label{eq-fixC}
\abs{C(k_m)}=\abs{C_{\bar{k}}^{F^m}}=1+q^m-a_m(C) \;\; \text{where} \;\; a_m(C) = \mathrm{Tr} \bigl(F^m,H^1(C_{\bar{k}},\Ql) \bigr).
\end{equation}
If we consider equations \eqref{eq-sigma} and \eqref{eq-equiv} in view of equation \eqref{eq-fixC} we find that
$$\abs{H_{g,n}^{F  \sigma}} = \sum_{[C] \in \Hh_{g}(k)/\cong_k} \frac{1}{\abs{\mathrm{Aut}_k(C)}} \, f_{\sigma}(q,a_1(C),\ldots,a_n(C)),$$ where $f_{\sigma}(x_0,\ldots,x_n)$ is a polynomial with coefficients in $\Z$. Give the variable $x_i$ degree $i$. Then there is a unique monomial in $f_{\sigma}$ of highest degree, namely $x_{1}^{\lambda_1} \cdots x_{\nu}^{\lambda_{\nu}}$. The numbers which we will pursue will be the following. 
\begin{dfn} \label{dfn-a} For $g \geq 2$ and any partition $\lambda$ define
\begin{equation} \label{eq-a}
a_{\lambda}|_g := \sum_{[C] \in \Hh_{g}(k)/\cong_k}  \frac{1}{\abs{\mathrm{Aut}_k(C)}} \, \prod_{i=1}^{\nu} a_i(C)^{\lambda_i}.
\end{equation} 
This expression will be said to have \emph{weight} $\abs{\lambda}$. Let us also define 
$$a_0|_g := \sum_{[C] \in \Hh_g(k)/\cong_k} \frac{1}{\abs{\mathrm{Aut}_k(C)}},$$ an expression of weight $0$.
\end{dfn}

\section{Representatives of hyperelliptic curves in odd characteristic} \label{sec-repr}
Assume that the finite field $k$ has an odd number of elements. The hyperelliptic curves of genus $g \geq 2$ are the ones endowed with a degree $2$ morphism to $\Pb^1$. This morphism induces a degree $2$ extension of the function field of $\Pb^1$. If we consider hyperelliptic curves defined over the finite field $k$ and choose an affine coordinate $x$ on $\Pb^1$, then we can write this extension in the form $y^2=f(x)$, where $f$ is a square-free polynomial with coefficients in $k$ of degree $2g+1$ or $2g+2$. At infinity, we can describe the curve given by the polynomial $f$ in the coordinate $t=1/x$ by $y^2=t^{2g+2}\, f(1/t)$. We will therefore let $f(\infty)$, which corresponds to $t=0$, be the coefficient of $f$ of degree $2g+2$.

\begin{dfn} \label{dfn-Pg} Let $P_g$ denote the set of square-free polynomials with coefficients in $k$ and of degree $2g+1$ or $2g+2$, and let $P'_g \subset P_g$ consist of the monic polynomials. Write $C_f$ for the curve corresponding to the element $f$ in $P_g$. 
\end{dfn}

By construction, there exists for each $k$-isomorphism class of objects in $\Hh_g(k)$ an $f$ in $P_g$ such that $C_f$ is a representative. Moreover, the $k$-isomorphisms between curves corresponding to elements of $P_g$ are given by $k$-isomorphisms of their function fields. By the uniqueness of the linear system $g^1_2$ on a hyperelliptic curve, these isomorphisms must respect the inclusion of the function field of $\Pb^1$. The $k$-isomorphisms are therefore precisely (see \cite[p. 126]{G-euler}) the ones induced by elements of the group $G:=\mathrm{GL}^{\mathrm{op}}_2(k) \times k^*/D$ where
$$D:=\{(\Bigl(\begin{array}{cc} a & 0 \\ 0 & a \end{array} \Bigr),a^{g+1}) : a \in k^*  \} \subset \mathrm{GL}^{\mathrm{op}}_2(k) \times k^*$$
and where an element 
$$\gamma= [(\Bigl( \begin{array}{cc} a & b \\ c & d \end{array} \Bigr),e)] \in G$$
induces the isomorphism
$$(x,y) \mapsto \left(\frac{ax+b}{cx+d},\frac{ey}{(cx+d)^{g+1}}\right).$$ 
This defines a left group action of $G$ on $P_g$, where $\gamma \in G$ takes $f \in P_g$ to $\tilde f \in P_g$, with
\begin{equation}\label{eq-lintrans} \tilde f(x)=\frac{(cx+d)^{2g+2}}{e^2} \, f \Bigl(\frac{ax+b}{cx+d}\Bigr).
\end{equation}

\begin{ntn} Let us put $I:=1/\abs{G}=(q^3-q)^{-1}(q-1)^{-1}$.
\end{ntn}

\begin{dfn}Let $\chi_{2,m}$ be the quadratic character on $k_m$. Recall that it is the function that takes $\alpha \in k_m$ to $1$ if it is a square, to $-1$ if it is a nonsquare and to $0$ if it is $0$.  With a square or a nonsquare we will always mean a nonzero element.
\end{dfn}

\begin{lem} \label{lem-am} If $C_f$ is the hyperelliptic curve corresponding to $f \in P_g$ then
$$a_m(C_f) = -\sum_{\alpha \in \Pb^1(k_m)} \chi_{2,m} \bigl( f(\alpha) \bigr).$$
\end{lem}
\begin{pf} The fiber of $C_f \to \Pb^1$ over $\alpha \in  \Ab^1(k_m)$ will consist of two points defined over $k_m$ if $f(\alpha)$ is a square in $k_m$, no point if $f(\alpha)$ is a nonsquare in $k_m$, and one point if $f(\alpha)=0$. By the above description of $f$ in terms of the coordinate $t=1/x$, the same holds for $\alpha=\infty$. The lemma now follows from equation~\eqref{eq-fixC}. 
\end{pf}

We will now rephrase equation \eqref{eq-a} in terms of the elements of $P_g$. By what was said above, the stabilizer of an element $f$ in $P_g$ under the action of $G$ is equal to $\Aut_k(C_f)$ and hence
\begin{multline} \label{eq-afp}
a_{\lambda}|_g = \sum_{[f] \in P_g/G} \frac{1}{\abs{\mathrm{Stab}_G(f)}} \, \prod_{i=1}^{\nu} a_{i}(C_f)^{\lambda_i} =\\  =  \frac{1}{\abs{G}}  \sum_{f \in P_g} \prod_{i=1}^{\nu} a_{i}(C_f)^{\lambda_i} =  I  \sum_{f \in P_g} \prod_{i=1}^{\nu} \Bigl(-\sum_{\alpha \in \Pb^1(k_{i})} \chi_{2,i} \bigl( f(\alpha) \bigr) \Bigr)^{\lambda_i}.
\end{multline} 
This can up to sign be rewritten as
\begin{equation} \label{eq-af}
I \sum_{f \in P_g} \sum_{(\alpha_{1,1}, \ldots,\alpha_{\nu,\lambda_{\nu}}) \in S} \prod_{i=1}^{\nu} \prod_{j=1}^{\lambda_i} \chi_{2,i} \bigl( f(\alpha_{i,j}) \bigr),
\end{equation}
where $S:=\prod_{i=1}^{\nu} \Pb^1(k_i)^{\lambda_i}$, in other words, $\alpha_{i,j} \in \Pb^1(k_{i})$ for each $1 \leq i \leq \nu$ and $1 \leq j \leq \lambda_i$. The sum \eqref{eq-af} will be split into parts for which we, in Section~\ref{sec-u}, will find recurrence relations in $g$.

\begin{dfn} For any tuple $\nn=(n_1, \ldots, n_m) \in \N_{\geq 1}^m$, let the set $A(\nn)$ consist of the tuples $\alpha=(\alpha_1, \ldots, \alpha_m) \in \prod_{i=1}^m \Pb^1(k_{n_i})$ such that for any $1\leq i,j \leq m$ and any $s \geq 0$, 
$$ F^s(\alpha_i)=\alpha_j \implies n_i|s \;  \text{and} \; i=j. $$
Let us also define $A'(\nn):=A(\nn) \cap \prod_{i=1}^m \Ab^1(k_{n_i})$.
\end{dfn}

\begin{dfn} Let $\mathcal{N}_m$ denote the set of pairs $(\nn;\rr)$ such that $\nn=(n_1,\ldots,n_m) \in \N_{\geq 1}^m$ and $\rr=(r_1,\ldots, r_m) \in \{1,2\}^m$.
\end{dfn}

\begin{dfn} \label{def-ug} For any $g \geq -1$, $(\nn;\rr) \in \mathcal{N}_m$ and $\alpha=(\alpha_1,\ldots,\alpha_m)\in A(\nn)$ define
$$ u_{g,\alpha}^{(\nn;\rr)} :=  I \sum_{f \in P_g} \prod_{i=1}^{m} \chi_{2,n_i} \bigl( f(\alpha_{i}) \bigr)^{r_i}$$
and 
$$u_{g}^{(\nn;\rr)} :=  \sum_{\alpha \in A(\nn)} u_{g,\alpha}^{(\nn;\rr)}.$$
\end{dfn}

\begin{constrlem} \label{lem-decomp}
For each $\lambda$, there are positive integers $c_1,\ldots,c_s$ and $m_1,\ldots,m_s$, and moreover pairs $(\nn^{(i)};\rr^{(i)}) \in \mathcal{N}_{m_i}$ for each $1 \leq i \leq s$, such that for any finite field $k$,
$$a_{\lambda}|_g=\sum_{i=1}^s c_i \, u_g^{(\nn^{(i)};\rr^{(i)})}.$$
\end{constrlem}
\begin{pf} The lemma will be proved by writing the set $S$ as a disjoint union of parts that only depend upon the partition $\lambda$, and which therefore are independent of the chosen finite field $k$.

For each positive integer $i$, let $i=d_{i,1} > \ldots > d_{i,\delta_i}=1$ be the divisors of~$i$.
\begin{itemize}
\item [$\star$] For each $1 \leq i \leq \nu$, let $T_{i,1},\ldots,T_{i,\delta_i}$ be an ordered partition of the set $\{1,\ldots,\lambda_i\}$ into (possibly empty) subsets.
\item [$\star$] For each $1 \leq i \leq \nu$ and each $1 \leq j \leq \delta_i$, let $Q_{i,j,1},\ldots, Q_{i,j,\kappa_{i,j}}$ be an unordered partition (where $\kappa_{i,j}$ is arbitrary) of the set $T_{i,j}$ into non-empty subsets.
\end{itemize}
From such a choice of partitions we define a subset $S'=S({\{T_{i,j} \},\{Q_{i,j,k} \}})$ of $S$ consisting of the tuples $(\alpha_{1,1},\ldots,\alpha_{\nu,\lambda_{\nu}}) \in S$ fulfilling the following two properties.
\begin{itemize}
\item [$\star$] If $x \in T_{i,j}$ then: $\alpha_{i,x} \in k_j \; \text{and} \; \forall s<j,  \alpha_{i,x} \notin k_s$. 
\item [$\star$] If $x \in Q_{i,j,k}$ and $x' \in Q_{i',j',k'}$ then:
$$\exists s : F^s(\alpha_{i,x})=\alpha_{i',x'}\iff  (i,j,k)=(i',j',k').$$
\end{itemize}

Define $\nn$ to be equal to the tuple
$$(\overbrace{d_{1,1},\ldots,d_{1,1}}^{\kappa_{1,1}},\overbrace{d_{1,2},\ldots,d_{1,2}}^{\kappa_{1,2}},\ldots,\overbrace{d_{1,\delta_1},\ldots,d_{1,\delta_1}}^{\kappa_{1,\delta_1}},\overbrace{d_{2,1},\ldots,d_{2,1}}^{\kappa_{2,1}},\ldots,\overbrace{d_{\nu,\delta_{\nu}},\ldots,d_{\nu,\delta_{\nu}}}^{\kappa_{\nu,\delta_{\nu}}}).$$
Let $\rho_{i,j,k}$ be equal to $2$ if either $i/d_{i,j}$ or $\abs{Q_{i,j,k}}$ is even, and $1$ otherwise. Define $\rr$  to be equal to
$$(\rho_{1,1,1},\rho_{1,1,2},\ldots,\rho_{1,1,\kappa_{1,1}},\rho_{1,2,1},\ldots,\rho_{1,\delta_1,\kappa_{1,\delta_1}},\rho_{2,1,1},\ldots,\rho_{\nu,\delta_{\nu},\kappa_{\nu,\delta_{\nu}}}).$$

The equality 
$$u_g^{(\nn;\rr)}=I \sum_{f \in P_g} \sum_{(\alpha_{1,1}, \ldots,\alpha_{\nu,\lambda_{\nu}}) \in S'} \prod_{i=1}^{\nu} \prod_{j=1}^{\lambda_i} \chi_{2,i} \bigl( f(\alpha_{i,j}) \bigr)$$
is clear in view of the following three simple properties of the quadratic character.

\begin{itemize}
\item[$\star$] Say that $\alpha \in \Pb^1(k_{s})$, then if $\tilde s/s$ is even we have $\chi_{2,\tilde s} \bigl( f(\alpha) \bigr) = \chi_{2,s} \bigl( f(\alpha) \bigr)^2$ and if $\tilde s / s$ is odd we have $\chi_{2,\tilde s} \bigl( f(\alpha) \bigr) = \chi_{2,s} \bigl( f(\alpha) \bigr)$. 
\item[$\star$] If for any $\alpha, \beta \in \Pb^1$ we have $F^s(\alpha)=\beta$ for some $s$, then $\chi_{2,i} \bigl( f(\alpha) \bigr) = \chi_{2,i} \bigl( f(\beta) \bigr)$ for all $i$. 
\item[$\star$] Finally, for any $\alpha \in \Pb^1$ and any $s$, we have $\chi_{2,s} \bigl( f(\alpha) \bigr)^r = \chi_{2,s} \bigl( f(\alpha) \bigr)^2$ if $r$ is even and $\chi_{2,s} \bigl( f(\alpha) \bigr)^r = \chi_{2,s} \bigl( f(\alpha) \bigr)$ if $r$ is odd.
\end{itemize} 
The lemma now follows directly from the fact that the sets $S(\{T_{i,j} \},\{Q_{i,j,k} \}) \subset S$ (for different choices of partitions $\{T_{i,j} \}$ and $\{Q_{i,j,k}\}$) are disjoint and cover $S$. \end{pf}

The set of data $\{(c_i, (\nn^{(i)}; \rr^{(i)}))\}$ resulting from the procedure given in the proof of Construction-Lemma \ref{lem-decomp} is, after assuming the pairs $(\nn^{(i)}; \rr^{(i)})$ to be distinct, unique up to simultaneous reordering of the elements of $\nn^{(i)}$ and $\rr^{ (i)}$ for each $i$, and it will be called \emph{the decomposition of} $a_{\lambda}|_g$.

\begin{dfn} \label{dfn-gen} For a partition $\lambda$, the pair
$$(\nn;\rr)=\bigl((\overbrace{1,\ldots,1}^{\lambda_1},\overbrace{2,\ldots,2}^{\lambda_2},\ldots,\overbrace{\nu,\ldots,\nu}^{\lambda_{\nu}});(1,\ldots,1) \bigr)$$
will appear in the decomposition of $a_{\lambda}|_g$ (corresponding to the partitions $T_{i,1}=\{1,\ldots,\lambda_i\}$ for $1\leq i \leq \nu$, and $Q_{i,1,k}=\{k\}$ for $1\leq i,k \leq \nu$) with coefficient equal to $1$, and it will be called the \emph{general case}. All other pairs $(\nn;\rr)$ appearing in the decomposition of $a_{\lambda}|_g$ will be refered to as \emph{degenerations} of the general case.
\end{dfn}

\begin{dfn} \label{dfn-n} For any $(\nn;\rr) \in \mathcal{N}_m$, the number $\abs{\nn}:=\sum_{i=1}^m n_i$ will be called the \emph{degree} of $(\nn;\rr)$.
\end{dfn}

\begin{lem} \label{lem-general} The general case is the only case in the decomposition of $a_{\lambda}|_g$ which has degree equal to the weight of $a_{\lambda}|_g$. \end{lem}
\begin{pf} If $(\nn;\rr)$ appears in the decomposition of $a_{\lambda}|_g$ and is associated to the partitions $\{T_{i,j} \}$ and $\{Q_{i,j,k}\}$, then $\abs{\nn}=\sum_{i=1}^{\nu} \sum_{j=1}^{\delta_i}\kappa_{i,j} d_{i,j}$. Since $\lambda_i = \sum_{j=1}^{\delta_i} \kappa_{i,j}$ and $1 \leq d_{i,j} \leq i$, the equality $\abs{\lambda}=\abs{\nn}$ implies that $\kappa_{i,1}=\lambda_i$ and $\kappa_{i,j}=0$ if $j \neq 1$.
\end{pf}

\begin{lem} \label{lem-decr}  If $(\nn;\rr)$ appears in the decomposition of $a_{\lambda}|_g$ then $\sum_{i=1}^m r_i n_i \leq \abs{\lambda}$ and these two numbers have the same parity.  \end{lem}
\begin{pf} If $(\nn;\rr)$ appears in the decomposition of $a_{\lambda}|_g$ and is associated to the partitions $\{T_{i,j} \}$ and $\{Q_{i,j,k}\}$, then $\sum_{i=1}^{m} r_i n_i=\sum_{i=1}^{\nu} \sum_{j=1}^{\delta_i} \sum_{k=1}^{\kappa_{i,j}} \rho_{i,j,k} d_{i,j}$. 

Let us prove the lemma by induction on $m$, starting with the case that $m=\sum_{i=1}^{\nu} \lambda_i$. In this case we must have $\abs{Q_{i,j,k}}=1$ for all $1 \leq i \leq \nu$, $1\leq j \leq \delta_i$ and $1 \leq k \leq \kappa_{i,j}$, and hence $\rho_{i,j,k}$ is only equal to two if $i/d_{i,j}$ is even. This directly tells us that $\rho_{i,j,k} d_{i,j} \leq i$, and that these two  numbers have the same parity. Since $\lambda_i=\sum_{j=1}^{\delta_i} \kappa_{i,j}$, it follows that $\sum_{i=1}^{m} r_i n_i \leq \abs{\lambda}$ and that these two numbers have the same parity. 

Assume now that $m=k$ and that the lemma has been proved for all pairs $(\tilde{\nn};\tilde{\rr})$ with $\tilde{m} >k$. Since $m < \sum_{i=1}^{\nu} \lambda_i$ we know that there exists numbers $i_0,j_0,k_0$ such that $\abs{Q_{i_0,j_0,k_0}} \geq 2$. Let us fix an element $x \in Q_{i_0,j_0,k_0}$ and define a new pair $(\nn';\rr')$ associated to the partitions $\{T'_{i,j} \}$ and $\{Q'_{i,j,k}\}$ by putting:
\begin{itemize}
 \item[$\star$] $T'_{i,j}=T_{i,j}$ for all $1 \leq i \leq \nu$ and $1 \leq j \leq \delta_i$,
 \item[$\star$] $Q'_{i_0,j_0,k_0}=Q_{i_0,j_0,k_0} \setminus \{x\}$, 
 \item[$\star$] $\kappa'_{i_0,j_0}=\kappa_{i_0,j_0}+1$ and $Q'_{i_0,j_0,\kappa'_{i_0,j_0}}=\{x\}$, 
 \item[$\star$] $Q'_{i,j,k}=Q_{i,j,k}$ in all other cases.
\end{itemize}
The pair $(\nn',\rr')$ thus appears in the decomposition of $\lambda$, and $m'=k+1$. Moreover, we directly find that $\sum_{i=1}^{m} r_i n_i \leq \sum_{i=1}^{m'} r'_i n'_i$ and that these two numbers have the same parity. By the induction hypothesis the lemma is then also true for $(\nn;\rr)$. \end{pf}

\begin{exa} Let us decompose $a_{[2^2]}|_g$ starting with the general case:
\begin{multline*} a_{[2^2]}|_g = I \sum_{f \in P_g} \Bigl(-\sum_{\alpha \in \Pb^1(k_2)} \chi_{2,2}\bigl(f(\alpha) \bigr) \Bigr)^2 = I \sum_{f \in P_g} \sum_{\alpha,\beta \in \Pb^1(k_2)} \chi_{2,2}\bigl(f(\alpha) f(\beta) \bigr) =\\ = u_g^{((2,2);(1,1))}+2u_g^{((2,1);(1,2))}+2u_g^{((2);(2))}+u_g^{((1,1);(2,2))}+u_g^{((1);(2))}.
\end{multline*}
\end{exa}

\begin{exa} \label{exa-dec} The decomposition of $a_{[1^4,2]}|_g$, starting with the general case:
\begin{multline*}
a_{[1^4,2]}|_g=
-u_g^{((2,1,1,1,1);(1,1,1,1,1))} - 6u_g^{((2,1,1,1);(1,2,1,1))} - 3u_g^{((2,1,1);(1,2,2))} \\
- 4u_g^{((2,1,1);(1,1,1))} - u_g^{((2,1);(1,2))} - u_g^{((1,1,1,1,1);(2,1,1,1,1))} - 6u_g^{((1,1,1,1);(2,2,1,1))} \\
- 4u_g^{((1,1,1,1),(1,1,1,1))} - 3u_g^{((1,1,1);(2,2,2))} - 22u_g^{((1,1,1);(2,1,1))} \\
- 7u_g^{((1,1);(2,2))} - 8u_g^{((1,1);(1,1))} - u_g^{((1);(2))}.
\end{multline*}
\end{exa}

\subsection{The cases of genus $0$ and $1$} \label{sec-g01} We would like to have an equality of the same kind as in equation \eqref{eq-afp}, but for curves of genus $0$ and $1$. Every curve of genus $0$ or~$1$ has a morphism to $\Pb^1$ of degree $2$ and in the same way as for larger genera, it then follows that every $k$-isomorphism class of curves of genus $0$ or $1$ has a representative among the curves coming from polynomials in $P_0$ and $P_1$ respectively. But there is a difference, compared to the larger genera, in that for curves of genus $0$ or $1$ the $g^1_2$ is not unique. In fact, the group $G$ induces (in the same way as for $g \geq 2$) all $k$-isomorphisms between curves corresponding to elements of $P_0$ and $P_1$ that respect their given morphisms to $\Pb^1$ (i.e a fixed $g^1_2$), but not all $k$-isomorphisms between curves of genus $0$ or $1$ are of this form.

Let us, for all $r \geq 0$, define the category $\mathcal{A}_r$ consisting of tuples $(C,Q_0,\ldots,Q_r)$ where $C$ is a curve of genus $1$ defined over $k$ and the $Q_i$ are, not necessarily distinct, points on $C$ defined over $k$. The morphisms of $\mathcal{A}_r$ are, as expected, isomorphisms of the underlying curves that fix the marked points. Note that $\mathcal{A}_0$ is isomorphic to the category $\Mm_{1,1}(k)$. We also define, for all $r \geq 0$, the category $\mathcal{B}_r$ consisting of tuples $(C,L,Q_1,\ldots,Q_r)$ of the same kind as above, but where $L$ is a $g^1_2$. A morphism of $\mathcal{B}_r$ is an isomorphism $\phi$ of the underlying curves that fixes the marked points, and such that there is an isomorphism $\tau$ making the following diagram commute:
$$\begin{CD}
C  @>\phi>> C'\\
@VLVV @VVL'V\\
\Pb^1 @>\tau>> \Pb^1.
\end{CD}$$

Consider $P_1$ as a category where the morphisms are given by the elements of $G$. To every element of $P_1$ there corresponds, precisely as for $g \geq 2$, a curve $C_f$ together with a $g^1_2$ given by the morphism to $\Pb^1$, thus an element of $\mathcal{B}_0$. Since every morphism in $\mathcal{B}_0$ between objects corresponding to elements of $P_1$ is induced by an element of $G$, and since for every $k$-isomorphism class of an element in $\mathcal{B}_0$ there is a representative in $P_1$, the two categories $P_1$ and $\mathcal{B}_0$ are equivalent.

For all $r \geq 1$ there are equivalences of the categories $\mathcal{A}_r$ and $\mathcal{B}_r$ given by
$$ (C,Q_0,\ldots,Q_r) \mapsto (C,|Q_0+Q_1|,Q_1,\ldots,Q_r),$$ with inverse $$ (C,L,Q_1,\ldots,Q_r) \mapsto (C,|L-Q_1|,Q_1,\ldots,Q_r).$$
We therefore have the equality
$$\sum_{[X]\in \mathcal{A}_r/\cong_k}\frac{1}{\abs{\Aut_k(X)}}\, \prod_{i=1}^{\nu} a_{i}(C)^{\lambda_i} = \sum_{[Y]\in \mathcal{B}_r/\cong_k}\frac{1}{\abs{\Aut_k(Y)}}\, \prod_{i=1}^{\nu} a_{i}(C)^{\lambda_i}.$$

The Riemann hypothesis tells us that $\abs{a_r(C)} \leq 2g \sqrt{q^r}$, for any finite field $k$ with $q$ elements and for any curve $C$ defined over $k$ of genus~$g$. For genus $1$ this implies that $\abs{C(k)} \geq q+1-2 \sqrt{q} >0$, and thus every genus $1$ curve has a point defined over $k$. There is therefore a number $s$ such that  $1 \leq \abs{C(k)} \leq s$ for all genus $1$ curves $C$. As in the argument preceding equation \eqref{eq-equiv} we can take a representative $(C,Q_0,\ldots,Q_r)$ for each element of $\mathcal{A}_r/\cong_k$ and act with $\Aut_k(C,Q_0)$, respectively for each representative $(C,L,Q_1,\ldots,Q_r)$ of $\mathcal{B}_0/\cong_k$ act with $\Aut_k(C,L)$, and by considering the orbits and stabilizers we get
\begin{equation*} \sum_{j=1}^s j^r \sum_{\substack{[X]\in \mathcal{A}_0/\cong_k \\\abs{C(k)}=j }}\frac{1}{\abs{\Aut_k(X)}}\, \prod_{i=1}^{\nu} a_{i}(C)^{\lambda_i} = \sum_{j=1}^s j^r  \sum_{\substack{[Y]\in \mathcal{B}_0/\cong_k \\ \abs{C(k)}=j}}\frac{1}{\abs{\Aut_k(Y)}}\, \prod_{i=1}^{\nu} a_{i}(C)^{\lambda_i}.\end{equation*}
Since this holds for all $r \geq 1$ we can, by a Vandermonde argument, conclude that we have an equality as above for each fixed $j$. We can therefore extend Definition~\ref{dfn-a} to genus $1$ in the following way:
\begin{multline} \label{eq-g1a}
a_{\lambda}|_1:=\sum_{\substack{[(C,Q_0)] \in \\\Mm_{1,1}(k)/\cong_k}} \frac{1}{\abs{\mathrm{Aut}_k(C,Q_0)}} \, \prod_{i=1}^{\nu} a_{i}(C)^{\lambda_i} = \\ =\sum_{[f] \in P_1/G} \frac{1}{\abs{\mathrm{Stab}_G(f)}} \, \prod_{i=1}^{\nu} a_{i}(C_f)^{\lambda_i} = I  \sum_{f \in P_{1}} \prod_{i=1}^{\nu} a_{i}(C_f)^{\lambda_i},\end{multline}
which gives an agreement with equation \eqref{eq-afp}.

All curves of genus $0$ are isomorphic to $\Pb^1$ and $a_r(\Pb^1)=0$ for all $r \geq 1$. In this trivial case we just let equation~\eqref{eq-afp} be the definition of $a_{\lambda}|_0$.

\section{Recurrence relations for $u_g$ in odd characteristic} \label{sec-u}
This section will be devoted to finding, for a fixed finite field $k$ with an odd number of elements and for a fixed pair $(\nn;\rr) \in \mathcal{N}_m$, a recurrence relation for $u_g$. Notice that we will often suppress the pair $(\nn;\rr)$ in our notation and for instance write $u_g$ instead of $u_g^{(\nn;\rr)}$. 

Fix a nonsquare $t$ in $k$ and an $\alpha=(\alpha_1,\ldots,\alpha_m) \in A(\nn)$. Multiplying with the element $t$ gives a fixed point free action on the set $P_g$ and therefore
\begin{multline} \label{eq-parity}
u_{g,\alpha}  = I  \sum_{f \in P_g} \prod_{i=1}^m \chi_{2,n_i} \bigl( f(\alpha_{i})\bigr)^{r_{i}} =   I  \sum_{f \in P_g} \prod_{i=1}^m \chi_{2,n_i} \bigl( t \, f(\alpha_{i})\bigr)^{r_{i}} =\\ =   I  \sum_{f \in P_g} \prod_{i=1}^m  \chi_{2,n_i}(t)^{r_{i}} \, \chi_{2,n_i} \bigl(f(\alpha_{i}) \bigr)^{r_{i}} =   (-1)^{\sum_{i=1}^m r_i n_i} \, u_{g,\alpha}. \end{multline}
This computation and Lemmas \ref{lem-decomp} and \ref{lem-decr} proves the following lemma.

\begin{lem} \label{lem-odd} For any $g \geq -1$, $(\nn;\rr) \in \mathcal{N}_m$ and $\alpha \in A(\nn)$, if $\sum_{i=1}^m r_i n_i$ is odd then $u_{g,\alpha}=0$. Consequently, $a_{\lambda}|_g$ is equal to $0$ if it has odd weight. \end{lem}

Thus, the only interesting cases are those for which $\sum_{i=1}^m r_i n_i$ is even. 

\begin{rmk} The last statement of Lemma~\ref{lem-odd} can also be found as a consequence of the existence of the hyperelliptic involution. \end{rmk} 

We also see from equation~\eqref{eq-parity} that
\begin{equation} \label{eq-monic} u_{g,\alpha}  = I \, (q-1) \sum_{f \in P'_g} \prod_{i=1}^m \chi_{2,n_i} \bigl( f(\alpha_{i})\bigr)^{r_{i}} \;\;\; \text{if $\sum_{i=1}^m r_in_i$ is even}.\end{equation}

\begin{dfn} \label{def-Ug} Let $Q_g$ denote the set of all polynomials (that is, not necessarily square-free) with coefficients in $k$ and of degree $2g+1$ or $2g+2$, and let $Q'_g \subset Q_g$ consist of the monic polynomials. For a polynomial $h \in Q_g$ we let $h(\infty)$ be the coefficient of the term of degree $2g+2$ (which extends the earlier definition for elements in $P_g$). For any $g \geq -1$, $(\nn;\rr)\in \mathcal{N}_m$ and $\alpha \in A(\nn)$, define 
$$ U_{g,\alpha}^{(\nn;\rr)} := I  \sum_{h \in Q_g} \prod_{i=1}^{m} \chi_{2,n_i} \bigl( h(\alpha_{i}) \bigr)^{r_i},$$ 
$$U_g^{(\nn;\rr)} := \sum_{\alpha \in A(\nn)}  U_{g,\alpha}^{(\nn;\rr)} \quad \text{and} \quad \UU_g^{(\nn;\rr)} := \sum_{i=-1}^g U_{i}^{(\nn;\rr)}.$$
\end{dfn}

We will find an equation relating $U_g$ to $u_{i}$ for all $-1 \leq i \leq g$. Moreover, for $g$ large enough we will be able to compute $U_g$. Together, this will give us our recurrence relation for $u_g$.

With the same arguments as was used to prove equation \eqref{eq-monic} one shows that
\begin{equation} \label{eq-monicU} U_{g,\alpha}  = I \, (q-1)  \sum_{h \in Q'_g} \prod_{i=1}^m \chi_{2,n_i}  \bigl(h(\alpha_{i})\bigr)^{r_{i}} \;\;\; \text{if $\sum_{i=1}^m r_in_i$ is even}.\end{equation}

\begin{dfn} \label{dfn-b}
For any $\alpha=(\alpha_1,\ldots,\alpha_m) \in A'(\nn)$, let $b_j=b^{\nn}_j$ be the number of monic polynomials $l$ of degree $j$ such that $l(\alpha_i)$ is nonzero for all $i$. Let us also put $\bb_j=\bb_j^{\nn}:=\sum_{i=0}^j b^{\nn}_i$. 
\end{dfn}

\begin{lem} For each $j \geq 0$ and $\nn \in \N_{\geq 1}^m$, we have the equality
\begin{equation} \label{eq-b}
b_j=q^j+\sum_{i=1}^j (-1)^{i}  \sum_{\substack{1 \leq m_1 <
\ldots < m_i \leq m \\ \sum_{l=1}^i n_{m_l} \leq j}} q^{j-\sum_{l=1}^i
n_{m_l}} \end{equation} 
from which it follows that $b_j$ does not depend upon the choice of $\alpha \in A'(\nn)$.
\end{lem}
\begin{pf} The numbers $b_j$ can be computed by inclusion-exclusion, where the choice of $1 \leq m_1 < \ldots < m_i \leq m$ corresponds to demanding the polynomial to be $0$ in the points $\alpha_{m_1}, \ldots, \alpha_{m_i}$. \end{pf}

\begin{ntn} For any $\alpha \in A'(\nn)$, let $p_{\alpha_i}$ denote the minimal polynomial of $\alpha_i$ and put $p_{\alpha}:=\prod_{i=1}^m p_{\alpha_i}$.
\end{ntn}

\begin{lem} \label{lem-onetoone} For any $\alpha \in A'(\nn)$ there is a one-to-one correspondence between polynomials $f$ defined over $k$ with $\deg(f) \leq \abs{\nn}-1$, and tuples $(f(\alpha_1),\ldots,f(\alpha_m)) \in \prod_{i=1}^m k_{n_i}$.
\end{lem}
\begin{pf}  For any $\alpha \in A'(\nn)$ we have $\deg(p_{\alpha_i})=n_i$ and $\gcd(p_{\alpha_i},p_{\alpha_j})=1$ if $i \neq j$. The lemma now follows from the Chinese remainder theorem, which tells us that the morphism $k[x]/p_{\alpha} \rightarrow \prod_{i=1}^m k[x]/p_{\alpha_i} \cong \prod_{i=1}^m k_{n_i}$ given by $f(x) \mapsto (f(\alpha_1),\ldots,f(\alpha_m))$ is an isomorphism. \end{pf}

\begin{ntn} Let $R_j$ denote the set of polynomials of degree $j$ and let $R'_j$ be the subset containing the monic polynomials. \end{ntn}

We will divide into two cases.

\subsection{The case $\alpha \in A'(\nn)$} \label{sec-n2} Fix an element $\alpha \in A'(\nn)$. Any nonzero polynomial $h$ can be written uniquely in the form $h=f \, l^2$ where $f$ is a square-free polynomial and $l$ is a monic polynomial. This statement translates directly into the equality
\begin{equation*}\label{eq-Uua} U_{s,\alpha} = I  \sum_{j+k=s} \sum_{l \in R'_j} \sum_{f \in P_k} \prod_{i=1}^{m} \chi_{2,n_i} \bigl( f(\alpha_{i}) \bigr)^{r_i} \chi_{2,n_i} \bigl( l(\alpha_{i}) \bigr)^{2r_i}  =\sum_{j=0}^{s+1} b_j u_{s-j,\alpha},\end{equation*}
because for any $\beta \in \Ab^1(k_s)$, $\chi_{2,s}\bigl((fl^2)(\beta)\bigr)=\chi_{2,s}\bigl(f(\beta)\bigr)$ if $l(\beta) \neq 0$.
Summing this equality over all $s$ between $-1$ and $g$ gives
\begin{equation}\label{eq-UUu} \UU_{g,\alpha} = \sum_{j=0}^{g+1} \bb_j u_{g-j,\alpha}. \end{equation} 

If $r_i=2$ for all $i$, then it follows from equation \eqref{eq-monicU} that
\begin{equation*} U_{s,\alpha} =  I \,(q-1)  \sum_{h \in Q'_s} \prod_{i=1}^{m} \chi_{2,n_i} \bigl( h(\alpha_{i}) \bigr)^{2} =  I \, (q-1) (b_{2s+2}+b_{2s+1}).\end{equation*}
Summing this equality over all $s$ between $-1$ and $g$ gives
\begin{equation} \label{eq-UUr0} \UU_{g,\alpha} = I\,(q-1) \bb_{2g+2} \quad \text{for $g\geq -1$ if $\forall i: r_i=2$.}\end{equation} 

In $\UU_{g,\alpha}$ we are summing over all polynomials $h$ of degree less than or equal to $2g+2$, and every $h$ can uniquely be written on the form $h_1+p_{\alpha} h_2$, with $\deg h_1 \leq \abs{\nn}-1$ and $\deg h_2 \leq 2g+2-\abs{\nn}$. Hence if $2g+2 \geq \abs{\nn}-1$ we find that
\begin{equation*} \UU_{g,\alpha} = I \, q^{2g+3-\abs{\nn}}  \sum_{s=1}^{\abs{\nn}-1} \sum_{h_1 \in R_s} \prod_{i=1}^{m} \chi_{2,n_i} \bigl( h_1(\alpha_{i}) \bigr)^{r_i}.
\end{equation*}
Using Lemma~\ref{lem-onetoone} we can reformulate this equality as
\begin{equation*} \UU_{g,\alpha} = I \, q^{2g+3-\abs{\nn}} \sum_{(\beta_1,\ldots,\beta_m) \in \prod_{i=1}^m k_{n_i} } \prod_{i=1}^m\chi_{2,n_i} ( \beta_i )^{r_i}.
\end{equation*}
For any $j$, half of the nonzero elements in $k_j$ are squares and half are nonsquares, and thus if $r_i=1$ for some $i$, we can conclude from this equality that
\begin{equation} \label{eq-UUr1} \UU_{g,\alpha}=0  \quad \text{for $g \geq (\abs{\nn}-3)/2$ if $\exists i: r_i=1$.}\end{equation}

\subsection{The case $\alpha \in A(\nn) \setminus A'(\nn)$} \label{sec-n1} Fix an element $\alpha \in A(\nn)\setminus A'(\nn)$. We can assume that $\alpha_1=\infty$, and then $\tilde \alpha:=(\alpha_2,\ldots,\alpha_m) \in A'(\tilde \nn)$ where $\tilde \nn:=(n_2,\ldots,n_m)$. 

If $h \in Q_g$ and $f \in P_j$ such that $h=f\,l^2$ for some monic polynomial $l$ (which is then unique), then $h(\infty)=f(\infty)$, because the coefficient of $h$ of degree $2g+2$ must equal the coefficient of $f$ of degree $2j+2$. As in Section~\ref{sec-n2} we get
\begin{equation}\label{eq-InfUua} U_{g,\alpha} = I  \sum_{j+k=g} \sum_{l \in R'_j} \sum_{f \in P_k} f(\infty)\prod_{i=2}^{m} \chi_{2,n_i} \bigl( f(\alpha_{i}) \bigr)^{r_i} \chi_{2,n_i} \bigl( l(\alpha_{i}) \bigr)^{2r_i}  = \sum_{j=0}^{g+1} b_j^{\tilde \nn} u_{g-j,\alpha}.\end{equation}

If $\sum_{i=1}^m r_i n_i$ is even, equation \eqref{eq-monicU} and the definition of $h(\infty)$ shows that
\begin{equation} \label{eq-monr0}
U_{g,\alpha} =  I \, (q-1) \sum_{h \in R'_{2g+2}} \prod_{i=2}^{m} \chi_{2,n_i} \bigl( h(\alpha_{i}) \bigr)^{r_i}.
\end{equation}
If $r_i=2$ for all $i$, then equation~\eqref{eq-monr0} tells us that
\begin{equation} \label{eq-InfUr0} U_{g,\alpha}= I \,(q-1)  b^{\tilde \nn}_{2g+2} \quad \text{for $g\geq -1$, $\forall i: r_i=2$}.\end{equation}

If $2g+2 \geq \abs{\nn}-1$, an element $h \in R'_{2g+2}$ can be written uniquely as $h=h_1+p_{\tilde \alpha}h_2$, where $\deg(h_1) \leq \abs{\nn}-2$, $\deg(h_2) \geq 0$ and $h_2$ monic. In the same way as in Section~\ref{sec-n2} we can (if $\sum_{i=1}^m r_i n_i$ is even) use this together with equation~\eqref{eq-monr0} and Lemma~\ref{lem-onetoone} to conclude that
\begin{equation} \label{eq-InfUr1} U_{g,\alpha}=0  \quad \text{for $g \geq (\abs{\nn}-3)/2$, $\exists i: r_i=1$},\end{equation}
which of course also holds if $\sum_{i=1}^m r_i n_i$ is odd by Lemma~\ref{lem-odd} and equation \eqref{eq-InfUua}.

\begin{rmk} Fix an $\alpha \in A(\nn)$. If there is an element $\beta \in \Ab^1(k)$ such that $\beta \notin \{\alpha_1,\ldots,\alpha_n\}$, then  $T(\alpha):=(T(\alpha_1),\ldots, T(\alpha_n))$ is in $A'(\nn)$, where $T$ is the projective transformation of $\Pb^1_k$ defined by $x \mapsto \beta x/(x-\beta)$.

In the notation of equation~\eqref{eq-lintrans}, $\chi_{2,n_i}\bigl(f(T(\alpha_i))\bigr)=\chi_{2,n_i}\bigl(\tilde f(\alpha_i) \bigr)$ (with $e=1$). Since this induces a permutation of $P_g$, we find that $u_{g,\alpha}=u_{g,T(\alpha)}$ and similarily that $U_{g,\alpha}=U_{g,T(\alpha)}$. So, if $q \geq \abs{\nn}$, then equations \eqref{eq-UUu}, \eqref{eq-UUr0} and \eqref{eq-UUr1} will also hold for $\alpha \in A(\nn) \setminus A'(\nn)$. By Lemma~\ref{lem-b} in the next section, we will see that this is true even if $q < \abs{\nn}$.
\end{rmk}

\subsection{The two cases joined} \label{sec-join}
In this section we will put the results of the two previous sections together using the following lemma. 

\begin{lem} \label{lem-b} For any $\tilde \nn=(n_2,\ldots,n_m)$, if $\nn=(1,n_2,\ldots,n_m)$ then $\bb^{\nn}_j=b^{\tilde \nn}_j$.
\end{lem}
\begin{pf}
Fix any tuple $\nn=(n_1,\ldots,n_m)$ and put $n:=\abs{\nn}$. If we let $t_i=q^{n_i}$ in the formula
$$\prod_{i=1}^m(t_i-1)=t_1 \cdots t_m + \sum_{i=1}^m (-1)^{i}  \sum_{1 \leq m_1 < \ldots < m_i \leq m} t_1 \cdots t_m \, \frac{1}{t_{m_1}} \cdots \frac{1}{t_{m_i}},$$
then the right hand side is equal to the right hand side of equation \eqref{eq-b}, and hence
\begin{equation} \label{eq-bn} \prod_{i=1}^m(q^{n_i}-1)=b^{\nn}_{n}. \end{equation}

Say that $b^{\nn}_j=\sum_{i=0}^j c^{\nn}_{j,i} q^i$ and $\bb^{\nn}_j=\sum_{i=0}^j \hat{c}^{\nn}_{j,i} q^i$. If $i\leq j$ then equation~\eqref{eq-b} implies that $c^{\nn}_{j,i}=c^{\nn}_{n,n+i-j}$ and hence $\hat{c}^{\nn}_{j,i}=\sum_{s=0}^j c^{\nn}_{n,n+i-s}$. By equation \eqref{eq-bn} we know that $q-1$ divides $b^{\nn}_n$, and if $b^{\nn}_{n}/(q-1)=\sum_{i=0}^{n-1} d_i q^i$ then $\hat{c}^{\nn}_{j,i}=d_{n-1+i-j}$. 

So, if $n_1=1$ and $\tilde \nn=(n_2,\ldots,n_m)$ then $b^{\nn}_{n}/(q-1)=b^{\tilde \nn}_{n-1}$ and thus $\hat{c}^{\nn}_{j,i}=c^{\tilde \nn}_{n-1,n-1+i-j}=c^{\tilde \nn}_{j,i}$.\end{pf}

\begin{ntn} Let us write $J:= I \, (q-1) \, \abs{A(\nn)}$. \end{ntn}

\begin{thm} \label{thm-rec1} For any pair $(\nn;\rr) \in \mathcal{N}_m$,
\begin{equation*} \sum_{j=0}^{g+1} \bb_j u_{g-j} = \begin{cases} J \, \bb_{2g+2} & \text{if $\forall i : r_i=2 $, $g \geq -1;$} \\ 0 & \text{if $\exists i : r_i=1$, $g \geq \frac{\abs{\nn}-3}{2}$.} 
\end{cases}  \end{equation*} 
\end{thm}
\begin{pf} The theorem follows from combining equations \eqref{eq-UUu}, \eqref{eq-UUr0}, \eqref{eq-UUr1} and equations \eqref{eq-InfUua}, \eqref{eq-InfUr0}, \eqref{eq-InfUr1}, using Lemma~\ref{lem-b}. \end{pf}

Note that with this theorem we can, for any $(\nn;\rr) \in \mathcal{N}_m$ such that $r_i=2$ for all~$i$, compute $u_g$ for any $g$. Moreover, for any pair $(\nn;\rr)$ we can compute $u_g$ for any $g$, if we already know $u_g$ for all $g < (\abs{\nn}-3)/2$.

\begin{lem} \label{lem-bcoeff} For any $\nn$, $q-1$ divides $b^{\nn}_{\abs{\nn}}$, and if we write $b^{\nn}_{\abs{\nn}}/(q-1)=\sum_{i=0}^{\abs{\nn}-1} d_i q^i$ then $\bb_j-q \bb_{j-1}=d_{\abs{\nn}-1-j}$.
\end{lem}
\begin{pf} The first claim is shown in the proof of Lemma~\ref{lem-b}. Using the notation of that proof we find that $\bb_j-q\bb_{j-1}=
\sum_{i=0}^j d_{n-1+i-j}q^i - \sum_{i=0}^{j-1} d_{n+i-j} q^{i+1} =d_{n-1-j}$. Note that $d_{n-1-j}$ only depends upon $\nn$ and not on $q$.
\end{pf}

\begin{thm} \label{thm-rec2} For any pair $(\nn;\rr) \in \mathcal{N}_m$, 
\begin{equation*} \sum_{j=0}^{\min(|\nn|-1,g+1)} (\bb_j-q \bb_{j-1})
 u_{g-j} = \begin{cases} J \, (\bb_{2g+2}-q\bb_{2g}) & \text{if $\forall i : r_i=2$, $g \geq 0;$} \\ 0 & \text{if $\exists i : r_i=1$, $g \geq \frac{\abs{\nn}-1}{2}.$}\end{cases}
\end{equation*}
\end{thm}
\begin{pf} Let us temporarily put $F(s):= \sum_{j=0}^{s+1} \bb_j u_{s-j}$. From Lemma~\ref{lem-bcoeff} we find that $\bb_j-q \bb_{j-1}=0$ if $j > \abs{\nn}-1$. The theorem then follows from applying Theorem~\ref{thm-rec1} to the expression $F(g)-q F(g-1)$. \end{pf}

For $g \geq (\abs{\nn}-1)/2$, Theorem \ref{thm-rec2} presents us with a linear recurrence relation for $u_g$ which has coefficients that are independent of the finite field $k$. 

\begin{exa} \label{exa-u1} If $(\nn;\rr)=((2,1,1,1);(1,2,1,1))$ then $b^{\nn}_5/(q-1)=(q^2-1)(q-1)^2=q^4-2q^3+2q-1$. Applying Lemma~\ref{lem-bcoeff} and then Theorem \ref{thm-rec2} we get
$$u_g-2 u_{g-1}+2 u_{g-3} -u_{g-4}=0 \quad \text{for $g \geq 3$}.$$
\end{exa}

\begin{exa} \label{exa-u2} Let us compute $u_g$, for all $g \geq -1$, when $(\nn;\rr)=((1,1,1),(2,2,2))$. We have that $u_{-1} = J = 1$ and since $r_i=2$ for all $i$, Theorem~\ref{thm-rec2} gives the equality $u_0 = 2 u_{-1} + J(q^2-3q+1) = q^2-3q+3$. Applying Theorem~\ref{thm-rec2} again we get
$$u_g-2 u_{g-1} + u_{g-2} = q^{2g-1}(q-1)^3 \quad \text{for $g \geq 1$}.$$
Solving this recurrence relation gives
$$u_g^{((1,1,1);(2,2,2))}= \frac{q^{2g+3}(q-1)-(2g+2)(q^2-1)+3q+1}{(q+1)^2} \quad \text{for $g \geq -1$}.$$
\end{exa}

\section{Linear recurrence relations for $a_{\lambda}|_g$} \label{sec-reca}
\begin{rmk} \label{rmk-seqsum} From a sequence $v_n$ that fulfills a linear recurrence relation with characteristic polynomial $C$ we can, for any polynomial $D$, in the obvious way construct a linear recurrence relation for $v_n$ with characteristic polynomial $CD$. Thus, from two sequences $v_n$ and $w_n$ that each fulfill linear recurence relation with characteristic polynomial $C$ and $D$ respectively, we can construct a linear recurence relation for the sequence $v_n+w_n$ with characteristic polynomial $\mathrm{lcm}(C,D)$.
\end{rmk}

\begin{thm} \label{thm-chareq} By applying Theorem~\ref{thm-rec2} to each pair $(\nn;\rr)$ appearing in the decomposition (given by Lemma~\ref{lem-decomp}) of $a_{\lambda}|_g$, we get a linear recurrence relation for $a_{\lambda}|_g$. The characteristic polynomial $C(X)$ of this linear recurrence relation equals
\begin{equation} \label{eq-chareq} \frac{1}{X-1} \, \prod_{i=1}^{\nu}(X^{i}-1)^{\lambda_i}.\end{equation}
\end{thm}
\begin{pf} Fix any pair $(\nn;\rr)$ in the decomposition of $a_{\lambda}|_g$ and put $n=\abs{\nn}$. Lemma~\ref{lem-bcoeff} tells us that $\bb_j-q\bb_{j-1}$ is equal to the coefficient of $q^{n-1-j}$ in $b_n/(q-1)$. If $g \geq n-1$, then these numbers are also the coefficients in the recurrence relation given by Theorem~\ref{thm-rec2}. By equation~\eqref{eq-bn}, the characteristic polynomial $C_{(\nn;\rr)}$ of this linear recurrence relation is equal to $(\prod_{i=1}^{m}(X^{n_i}-1))/(X-1)$. 

We find that the linear recurrence relation in the general case (see Definition~\ref{dfn-gen}) will have characteristic polynomial equal to $C$. Moreover, we find (by their construction in the proof of Lemma~\ref{lem-decomp}) that if $(\nn;\rr)$ is a degenerate case then $C_{(\nn;\rr)} | C$. The theorem now follows from Remark~\ref{rmk-seqsum}. \end{pf}

Theorem~\ref{thm-chareq} tells us that if we can compute $a_{\lambda}|_g$ for $g < \abs{\lambda}-1$ then we can compute it for every $g$. But note that by considering the individual cases in the decomposition of $a_{\lambda}|_g$ we will do much better in Section~\ref{sec-results}, in the sense that we will be able to use information from curves of only genus $0$ and $1$ to compute $a_{\lambda}|_g$ for any $\lambda$ such that $\abs{\lambda} \leq 6$.

\begin{exa} For $\lambda=[1^4,2]$ the characteristic polynomial equals $(X-1)^4(X+1)$, so if $V_g$ is a particular solution to the linear recurrence relation for $a_{[1^4,2]}|_g$ then 
\begin{equation*} a_{[1^4,2]}|_g = V_g+A_3 g^3+A_2 g^2+A_1 g+A_0+B_0 (-1)^g, \end{equation*}
where $A_0$, $A_1$, $A_2$, $A_3$ and $B_0$ do not depend upon $g$.
\end{exa}

\section{Computing $u_0$} \label{sec-gzero}
In this section we will see that we can compute $u_0$ for any choice of a pair $(\nn;\rr) \in \mathcal{N}_m$. This is due to the fact that if $C$ is a curve of genus $0$ then, for all~$r$, $\abs{C(k_r)}=1+q^r$ or equivalently $a_r(C)=0$. 

\begin{constrlem} \label{lem-g0} For each $(\nn;\rr) \in \mathcal{N}_m$, there are numbers $c_1,\ldots,c_s$ and pairs $(\nn^{(1)};\rr^{(1)}),\ldots,(\nn^{(s)};\rr^{(s)})$, where $\rr^{(i)}=(2,\ldots,2)$ for all $i$, such that for any finite field $k$, 
$$u^{(\nn;\rr)}_0=\sum_{i=1}^s c_i u_0^{(\nn^{(i)};\rr^{(i)})}.$$ 
\end{constrlem}
\begin{pf} Fix a pair $(\nn;\rr)\in \mathcal{N}_m$. We will use induction over the number $n:=\abs{\nn}$, where the base case $n=0$ is trivial.

Let us put $(\tilde \nn;\tilde \rr)=((n_2,\ldots,n_m);(r_2,\ldots,r_m))$. For an $\tilde \alpha=(\alpha_2,\ldots,\alpha_m) \in A(\tilde \nn)$ let $\hat{\Pb}^1_{\tilde \alpha}(k_{i})$ be the set of all points in $\Pb^1(k_{i}) \setminus \{\alpha_2,\ldots,\alpha_m\}$ that are not defined over a proper subfield of $k_{i}$.
The set of $\alpha_1 \in \Pb^1(k_{n_1})$ such that $(\alpha_1,\ldots,\alpha_m) \in A(\nn)$ then equals
\begin{equation} \label{eq-setsg0}
\Pb^1(k_{n_1}) \setminus \bigl(\bigcup_{i | n_1} \hat{\Pb}_{\tilde \alpha}^1(k_{i})\bigcup_{n_i | n_1} \{\alpha_i,\ldots,F^{n_i-1} \alpha_i \} \bigr).
\end{equation}

Assume now that the lemma has been proved for all pairs of degree strictly less than $n$. By reordering the elements of the pair $(\nn;\rr)$ we can assume that $r_1=1$, because otherwise $\rr=(2,\ldots,2)$ and we are done. By applying equation~\eqref{eq-setsg0} we get
\begin{multline} \label{eq-g0}
I \sum_{\alpha \in A(\nn)} \prod_{i=1}^m \chi_{2,n_i} \bigl( f(\alpha_{i}) \bigr)^{r_{i}} = I  \sum_{\tilde \alpha\in A(\tilde \nn)}  \prod_{i=2}^{m} \chi_{2,n_i} \bigl( f(\alpha_{i}) \bigr)^{r_{i}} \cdot \\ \cdot \Bigl(-a_{n_1}(C_f) -\sum_{i |n_1} \sum_{\beta \in \hat{\Pb}_{\tilde \alpha}^1(k_i)} \chi_{2,n_1} \bigl( f(\beta) \bigr) - \sum_{n_i |n_1} n_i \, \chi_{2,n_1} \bigl( f(\alpha_{i}) \bigr) \Bigr).
\end{multline}
Let us put $(\nn^{(i)};\rr^{(i)})=((i,n_2,\ldots,n_m);(n_i/i,r_2,\ldots,r_m))$ for all $i$ that divides $n_i$ and $\tilde \rr^{(i)}=(r_2,\ldots,r_{i-1},r_i \, n_1/n_i,r_{i+1},\ldots,r_m)$ for all $n_i$ that divides $n_1$. Summing both sides of equation~\eqref{eq-g0} over polynomials $f \in P_0$ and using that $a_{n_1}(C_f)=0$ then gives
\begin{equation} \label{eq-decompg0} u^{(\nn;\rr)}_0 = -\sum_{i |n_1} u_0^{(\nn^{(i)};\rr^{(i)})}-\sum_{n_i |n_1} n_i  u_0^{(\tilde \nn; \tilde \rr^{(i)})}. \end{equation}
Since $\abs{\tilde \nn} <n$ and $\abs{\nn^{(i)}} <n$, the lemma follows by induction from equation~\eqref{eq-decompg0}. \end{pf}

\begin{exa} In the case $(\nn;\rr)=((6,6,3,1,1);(1,1,2,2,2))$, the first step in the procedure in the proof of Lemma~\ref{lem-g0} equals
\begin{multline*}
u_0^{(\nn;\rr)}=-u_0^{((6,3,3,1,1);(1,2,2,2,2))}-u_0^{((6,3,2,1,1);(1,2,1,2,2))}-u_0^{((6,3,1,1,1);(1,2,2,2,2))}\\ -5u_0^{((6,3,1,1);(1,2,2,2))}-6u_0^{((6,3,1,1);(2,2,2,2))}.
\end{multline*}
\end{exa}

\begin{exa} In the case $(\nn;\rr)=((4,1,1,1);(1,2,1,1))$, the procedure in the proof of Lemma~\ref{lem-g0} gives
\begin{multline*}
u_0^{(\nn;\rr)} =u_0^{((2,1,1);(2,2,2))}+u_0^{((1,1,1);(2,2,2))}+u_0^{((1,1);(2,2))}-u_0^{((2,1);(2,2))}-u_0^{((1);(2))} \end{multline*}
\end{exa}

\section{Results for weight up to $7$ in odd characteristic} \label{sec-results} 
We will in this section show that we, for any number $g$ and any finite field $k$ of odd characteristic, can compute all $a_{\lambda}|_g$ of weight at most~$7$. This is achieved by decomposing $a_{\lambda}|_g$ using Lemma~\ref{lem-decomp} and employing the recurrence relation of Theorem~\ref{thm-rec1} on the different parts. This involves finding the necessary base cases for the recurrence relations and that will be possible with the help of results on genus~$0$ curves obtained in Section~\ref{sec-gzero}, and on genus $1$ curves obtained in the article~\cite{Jonas2}.

We will write $a_{\lambda}|_{g,odd}$ and $u_{g,odd}$ to stress that all results are in the case of odd characteristic. See Section~\ref{sec-results2} for results in the case of even characteristic. 

\begin{exa} \label{exa-a0} Theorem~\ref{thm-rec1} is applicable even if the degree is $0$ (if considered as a case when $r_i=2$ for all $i$) and with $\bb_j=\sum_{i=0}^j q^i$. From Theorem~\ref{thm-rec1} we find that $a_0|_{0,odd}=Jq^2=q/(q^2-1)$ and again from Theorem~\ref{thm-rec1} that
$$a_0|_{g,odd}=J(q^{2g+2}-q^{2g})=q^{2g-1} \quad \text{for $g \geq 1$.}$$
This result can also be found in \cite[Proposition 7.1]{Brock}.
\end{exa}

\subsection{Degree at most $3$} \label{sec-dthree}
When the degree of the pair $(\nn;\rr)$ is at most $3$ we find using Theorem~\ref{thm-rec1} that we do not need any base cases to compute $u_{g}$ for every $g$.

\begin{exa} \label{exa-w2} Let us consider $(\nn;\rr)=((2);(1))$. We have $u_{-1}=J= 1/(q+1)$ and applying Theorem~\ref{thm-rec1} we get $u_{0} = -(q+1) u_{-1}= -1$. Theorem~\ref{thm-rec2} tells us that $u_g=-u_{g-1}$ for $g \geq 1$ and thus
$$u^{((2);(1))}_{g,odd} =(-1)^{g+1} \quad \text{for $g \geq 0$}.$$
\end{exa}

\begin{exa}The result for $a_{[2]}|_{g,odd}$ is
$$a_{[2]}|_{g,odd} = -u^{((2);(1))}_g - u^{((1);(2))}_g = (-1)^g - q^{2g} \quad \text{for $g \geq 0$}.$$
\end{exa}

\begin{exa}The result for $a_{[1^2]}|_{g,odd}$ is
$$a_{[1^2]}|_{g,odd} = u^{((1,1);(1,1))}_g + u^{((1);(2))}_g = -1 + q^{2g} \quad \text{for $g \geq 0$}.$$
\end{exa}

\begin{rmk} 
The result for $(q^2+1) \, a_0|_{g,odd}-a_{[2]}|_{g,odd}$ can be found in lecture notes by Bradley Brock and Andrew Granville from 28 July 2003.
\end{rmk}

\begin{exa} Consider the case $(\nn;\rr)=((1,1,1);(2,1,1))$. We have $u_{-1}=J=1$ and from Theorem~\ref{thm-rec1} we get $u_0=-(q-2)u_{-1}=-q+2$. Theorem~\ref{thm-rec2} gives the recurrence relation $u_{g}=2u_{g-1}-u_{g-2}$ for $g \geq 1$ and hence
$$u_{g,odd}^{((1,1,1);(2,1,1))}=g(-q+1)-q+2.$$
\end{exa}

\subsection{Degree $4$ or $5$}  \label{sec-dfive}
From Theorem~\ref{thm-rec1} we find that when the degree of the pair $(\nn;\rr)$ is $4$ or $5$ we need the base case of genus $0$. But the genus $0$ case is always computable using Lemma~\ref{lem-g0} and then Theorem~\ref{thm-rec1}, and hence the same is true for $u_g$ for all $g$.

\begin{exa} For $(\nn;\rr)=((2,1,1);(1,1,1))$ we have $u_{-1}=q$ and from Lemma~\ref{lem-g0} it follows that $$u_0^{((2,1,1);(1,1,1))}=-u_0^{((2,1);(1,2))}=u_0^{((1,1);(2,2))}+u_0^{((1);(2))}=q.$$ Using Theorem~\ref{thm-rec1} we get $u_1=-(q-1)u_0-(q^2-q-1)u_{-1}=-q^3+2q$. Solving the recurrence relation $u_{g}=u_{g-1}-u_{g-2}-u_{g-3}$ for $g \geq 2$, coming from Theorem~\ref{thm-rec2}, gives
$$u_{g,odd}^{((2,1,1);(1,1,1))}=\frac{1}{4} \, (q^3-q)(-2g+(-1)^g-1)+ q.$$
\end{exa}

\begin{exa} The result for $a_{[1^2,2]}|_{g,odd}$ is
\begin{multline*}
a_{[1^2,2]}|_{g,odd} = - u_g^{((2,1,1);(1,1,1))}-u_g^{((2,1);(1,2))}-u_g^{((1,1,1);(2,1,1))}\\- u_g^{((1,1);(2,2))} - 2u_g^{((1,1);(1,1))}-u_g^{((1);(2))} =\\=-\frac{q^{2g+2}-1}{q+1} -q^{2g} +\frac{1}{2}\, g(q^3+q-2) +\frac{1}{2} \begin{cases} 2q & \text{if $g \equiv 0$ mod $2$}\\q^3-q-2 & \text{if $g \equiv 1$ mod $2$}\end{cases}
\end{multline*}
\end{exa}

\subsection{Weight $6$}  \label{sec-wsix}
We will not be able to compute $u_{g}$ for all pairs $(\nn;\rr)$ of degree~$6$. But we will be able to compute $u_g$ for all pairs $(\nn;\rr)$ that are general cases in the decomposition of $a_{\lambda}|_g$ for $\lambda$'s of weight $6$. This will be sufficient to compute all $a_{\lambda}|_g$ of weight $6$, because we saw in Lemma~\ref{lem-general} that only the general case will have degree $6$ and therefore all degenerate cases are covered in Sections \ref{sec-dthree} and \ref{sec-dfive}.

Let $u_g$ be the general case in the decomposition of $a_{\lambda}|_g$. When the degree is equal to $6$ we see from Theorem~\ref{thm-rec1} that we need the base cases of genus $0$ and $1$ to compute $u_g$ for all $g$. As we know, we can always compute $u_0$ using Lemma~\ref{lem-g0}. For genus $1$, the numbers $a_{\lambda}|_1$ have been computed for weight up to $6$ by the author. This was done by embedding every genus $1$ curve with a given point as a plane cubic curve, see \cite[Section 15]{Jonas2}. Since we know all the degenerate cases in the decomposition of $a_{\lambda}|_1$ we can then compute the general case $u_1$.

\begin{exa} Let us deal with $(\nn;\rr)=((6);(1))$ which is the generic case in the decomposition of $a_{[6]}|_{g,odd}$ and for which we have $u_{-1}=J=q^3+q-1$. Using Lemma~\ref{lem-g0} we get
$$u_0^{((6);(1))}=-u_0^{((3);(2))}-u_0^{((2);(1))}-u_0^{((1);(2))}=-u_0^{((3);(2))}=-q^2.$$ Using the results of \cite[Section 15]{Jonas2} we find that $a_{[6]}|_1=q-1$. Decomposing $a_{[6]}|_g$ gives $a_{[6]}|_1=-u_1^{((6);(1))}-u_1^{((3);(2))}-u_1^{((2);(1))}-u_1^{((1);(2))}$. Thus, using Example \ref{exa-w2}, we get $u_1=-(q-1)-(q^4-q^2-q-1)-1-q^2=-q^4+1$. We can now apply Theorem~\ref{thm-rec1} which gives $u_2=-(q+1)u_1-(q^2+q+1)u_0-(q^3+q^2+q+1)u_{-1}=-q^6+q^2-q$, $u_3=-u_2-u_1-u_0-u_{-1}=q^6+q^4-q^3$ and $u_4=-u_3-u_2-u_1-u_0-u_{-1}=0$. If we then multiply the characteristic polynomial for the linear recurrence relation of $u_g$ by $X-1$ we get $u_{g}=u_{g-6}$ for all $g \geq 5$.
\end{exa}

\begin{exa}The result for $a_{[6]}|_{g,odd}$ is
\begin{multline*}
a_{[6]}|_{g,odd} =-u^{((6);(1))}_g-u^{((3);(2))}_g-u^{((2);(1))}_g-u^{((1);(2))}_g= -q^{2g}-\frac{q^{2g+3}(q-1)}{q^2-q+1}+\\+\frac{1}{q^2-q+1}  \begin{cases}q^2 & \text{if $g \equiv 0$ mod $3$}\\-q^2-1 & \text{if $g \equiv 1$ mod $3$}\\1 & \text{if $g \equiv 2$ mod $3$}\end{cases} +\begin{cases} q^2+1 & \text{if $g \equiv 0$ mod $6$}\\q^4-2 & \text{if $g \equiv 1$ mod $6$} \\ q^6-q^2+q+1 & \text{if $g \equiv 2$ mod $6$}
\\-q^6-q^4+q^3-1 & \text{if $g \equiv 3$ mod $6$} \\1 & \text{if $g \equiv 4$ mod $6$} \\-q^3-q & \text{if $g \equiv 5$ mod $6$} \end{cases}
\end{multline*}
\end{exa}

\begin{rmk} \label{rmk-polynomial} For any choice of $\lambda$ and $g$, consider $a_{\lambda}|_{g,odd}$ as a function of the number $q$ of elements of the finite field $k$ of odd characteristic. If $\lambda$ is of weight at most $7$ it follows from our computations that this function is a polynomial in the variable $q$.

This will not continue to hold when considering for instance $a_{[1^6]}|_3$, that is, also including finite fields of \emph{even} characteristic, see Example~\ref{exa-a16even}. But it will also not hold for instance for $a_{[1^{10}]}|_{1,odd}$, which for prime fields will be a polynomial function minus the Ramanujan $\tau$-function, compare \cite[Corollary 5.4]{G-res}. 
\end{rmk}

\section{Representatives of hyperelliptic curves in even characteristic} \label{sec-repreven}
Let $k$ be a finite field with an even number of elements. We will again describe the hyperelliptic curves of genus $g \geq 2$ defined over $k$ by their degree $2$ morphism to $\Pb^1$. If we choose an affine coordinate $x$ on $\Pb^1$ we can write the induced degree~$2$ extension of the function field of $\Pb^1$ in the form $y^2+h(x)y+f(x)=0$, where $h$ and $f$ are polynomials defined over $k$ that fulfill the following conditions:
\begin{gather} \label{eq-deg} 2g+1 \leq \max\bigl(2 \deg(h),\deg(f)\bigr) \leq 2g+2; \\ \label{eq-nons} \gcd(h,f'^2+fh'^2) = 1; \\ \label{eq-inf} t \nmid \gcd(h_{\infty},f_{\infty}'^2+f_{\infty}h_{\infty}'^2). 
\end{gather}
The last condition comes from the nonsingularity of the point(s) in infinity, around which the curve can be described in the variable $t=1/x$ as $y^2+h_{\infty}(t)y+f_{\infty}(t)=0$, where $h_{\infty}:=t^{g+1}h(1/t)$ and $f_{\infty}:=t^{2g+2}f(1/t)$. We therefore define $h(\infty)$ and $f(\infty)$ to be equal to the degree $g+1$ and $2g+2$ coefficient respectively. For a reference see for instance \cite[p. 294]{Liu}.

\begin{dfn} 
Let $P_g$ denote the set of pairs $(h,f)$ of polynomials defined over $k$, where $h$ is nonzero, that fulfill all three conditions \eqref{eq-deg}, \eqref{eq-nons} and \eqref{eq-inf}. Write $C_{(h,f)}$ for the curve corresponding to the element $(h,f)$ in $P_g$. 
\end{dfn}

To each $k$-isomorphism class of objects in $\Hh_g(k)$ there is a pair $(h,f)$ in $P_g$ such that $C_{(h,f)}$ is a representative. All $k$-isomorphisms between the curves represented by elements of $P_g$ are given by $k$-isomorphisms of their function fields, and since the $g^1_2$ of a hyperelliptic curve is unique the $k$-isomorphisms must respect the inclusion of the function field of $\Pb^1$.

Identify the set of polynomials $l(x)$ defined over $k$ and of degree at most $g+1$ with $k^{g+2}$, and define the group homomorphism
$$\phi_g:\mathrm{GL}^{\mathrm{op}}_2(k) \times k^* \rightarrow \Aut(k^{g+2}), \; \phi_g(\Bigl(\begin{array}{cc} a & b \\ c & d \end{array} \Bigr),e)\bigl(l(x)\bigr):=e^{-1}(cx+d)^{g+1} l\Bigl(\frac{ax+b}{cx+d}\Bigr).$$
Now define the group $G_g:=\bigl(k^{g+2}\rtimes_{\phi_g}(\mathrm{GL}^{\mathrm{op}}_2(k) \times k^*) \bigr)/D$ where
$$D:=\{(0,\Bigl(\begin{array}{cc} a & 0 \\ 0 & a \end{array} \Bigr),a^{g+1}) : a \in k^*  \} \subset k^{g+2}\rtimes_{\phi_g}(\mathrm{GL}^{\mathrm{op}}_2(k) \times k^*).$$
The $k$-isomorphisms between curves corresponding to elements of $P_g$ are then precisely the ones induced by elements of $G_g$ by letting
$$\gamma= [(l(x),\Bigl( \begin{array}{cc} a & b \\ c & d \end{array} \Bigr),e)] \in G_g $$
induce the isomorphism
$$(x,y) \mapsto \left(\frac{ax+b}{cx+d},\frac{e\bigl(y+l(x)\bigr)}{(cx+d)^{g+1}}\right).$$ 
This defines a left group action of $G_g$ on $P_g$, where $\gamma=[(l,\Lambda,e)] \in G_g$ takes $(h,f) \in P_g$ to $(\tilde h,\tilde f) \in P_g$, with
$$(\tilde h,\tilde f)= (\phi_g(\Lambda,e)(h),e^{-1} \phi_{2g}(\Lambda,e)(f)+ l\, \phi_g(\Lambda,e)(h)+l^2).$$

\begin{dfn} Let $\tau_m$ be the function that takes $(a,b) \in k_m^2$ to $1$ if the equation $y^2+ay+b$ has two roots defined over $k_m$, $0$ if it has one root and $-1$ if it has none. 
\end{dfn}

\begin{lem} If $C_{(h,f)}$ is the hyperelliptic curve corresponding to $(h,f) \in P_g$ then
\begin{equation*} a_m(C_{(h,f)})=-\sum_{\alpha \in \Pb^1(k_m)} \tau_m \bigl(h(\alpha),f(\alpha) \bigr).\end{equation*}
\end{lem}
\begin{pf} Follows in the same way as Lemma~\ref{lem-am}. \end{pf}

\begin{ntn} Let us put $I_g:=1/\abs{G_g}=q^{-(g+2)}(q^3-q)^{-1}(q-1)^{-1}$. \end{ntn}

In the same way as in the case of odd characteristic we get the equality
\begin{equation*} a_{\lambda}|_g = I_g  \sum_{(h,f) \in P_g} \prod_{i=1}^{\nu} \Bigl(-\sum_{\alpha \in \Pb^1(k_{i})} \tau_{i} \bigl(h(\alpha),f(\alpha) \bigr) \Bigr)^{\lambda_i}.\end{equation*}

All results of Section~\ref{sec-g01} are independent of the characteristic and hence we extend the definition of $a_{\lambda}|_g$ to genus $0$ and $1$ in the same way as in that section.

\begin{dfn} \label{def-ugeven} 
For any $g \geq -1$, $(\nn;\rr) \in \mathcal{N}_m$ and $\alpha \in A(\nn)$ define 
$$u_{g,\alpha}^{(\nn;\rr)} :=  I_g  \sum_{(h,f) \in P_g} \prod_{i=1}^m \tau_{n_i} \bigl(h(\alpha_i),f(\alpha_{i}) \bigr)^{r_{i}} $$
and 
$$u_{g}^{(\nn;\rr)}:= \sum_{\alpha \in A(\nn)} u_{g,\alpha}^{(\nn;\rr)}. $$
\end{dfn}

\begin{constrlem} \label{lem-decompeven} For each $\lambda$ we have (in even characteristic) the same decomposition of $a_{\lambda}|_g$ as given by Construction-Lemma~\ref{lem-decomp}.
\end{constrlem}
\begin{pf} The following properties of $\tau_m$ for $(h,f) \in P_g$ correspond precisely to the ones for the quadratic character.
\begin{itemize}
\item[$\star$] Say that $\alpha \in \Pb^1(k_{s})$, then $\tau_{\tilde s} \bigl(h(\alpha),f(\alpha) \bigr) = \tau_{s} \bigl(h(\alpha), f(\alpha) \bigr)^2$ if $\tilde s/s$ is even, and $\tau_{\tilde s} \bigl(h(\alpha), f(\alpha) \bigr) = \tau_{s} \bigl(h(\alpha), f(\alpha) \bigr)$ if $\tilde s / s$ is odd. 
\item[$\star$] If for any $\alpha, \beta \in \Pb^1$ we have $F^s(\alpha)=\beta$ for some $s$, then $\tau_{i} \bigl(h(\alpha),  f(\alpha) \bigr) = \tau_{i} \bigl(h(\alpha),  f(\beta) \bigr)$ for all $i$. 
\item[$\star$] Finally, for any $\alpha \in \Pb^1$ and any $s$, $\tau_{,s} \bigl(h(\alpha),  f(\alpha) \bigr)^r = \tau_{s} \bigl(h(\alpha),  f(\alpha) \bigr)^2$ if $r$ is even and $\tau_{s} \bigl(h(\alpha),  f(\alpha) \bigr)^r = \tau_{s} \bigl(h(\alpha),  f(\alpha) \bigr)$ if $r$ is odd.
\end{itemize} 
With this established we can use the same proof as for Construction-Lemma~\ref{lem-decomp}.
\end{pf}

Since the decompositions are the same, Lemmas~\ref{lem-general} and \ref{lem-decr} also hold in even characteristic.

\section{Recurrence relations for $u_g$ in even characteristic} \label{sec-ueven}
Analogously to Section~\ref{sec-u}, this section will be devoted to finding for a fixed pair $(\nn;\rr) \in \mathcal{N}_m$, a recurrence relation for $u_g$. Fix an $s \in k$ which does not lie in the set $\{r^2+r: r \in k \}$, that is, such that $\tau_1(1,s)=-1$. We define an involution on $P_g$ sending $(h,f)$ to $(h,f+s \, h^2)$. This involution is fixed point free and hence
\begin{multline*} 
u_{g,\alpha}  = I_g  \sum_{(h,f) \in P_g} \prod_{i=1}^m \tau_{n_i} \bigl( h(\alpha_i), f(\alpha_{i})\bigr)^{r_{i}} = \\ = I_g  \sum_{(h,f) \in P_g} \prod_{i=1}^m \tau_{n_i} \bigl( h(\alpha_i), f(\alpha_{i}) + s \, h^2(\alpha_i)\bigr)^{r_{i}} =(-1)^{\sum_{i=1}^m r_i n_i} u_{g,\alpha}. \end{multline*}
Thus, Lemma~\ref{lem-odd} also holds in the case of even characteristic.

\begin{dfn} Let $Q_g$ denote the set of pairs $(h,f)$ of polynomials over $k$, where $h$ is nonzero and $h$, $f$ are of degree at most $g+1$, $2g+2$ respectively. Extending the definition for $P_g$ above to a pair $(h,f) \in Q_g$, let $h(\infty)$ and $f(\infty)$ be equal to the degree $g+1$ and $2g+2$ coefficient of $h$ and $f$ respectively. For any $g \geq -1$, $(\nn;\rr) \in \mathcal{N}_m$ and $\alpha \in A(\nn)$ define
$$\UU_{g,\alpha}^{(\nn;\rr)} := I_g   \sum_{(h,f) \in Q_g} \prod_{i=1}^m \tau_{n_i} \bigl(h(\alpha_i),f(\alpha_{i}) \bigr)^{r_{i}}$$ 
and 
$$\UU_{g}^{(\nn;\rr)} :=\sum_{\alpha \in A(\nn)}\UU_{g,\alpha}^{(\nn;\rr)}.$$
\end{dfn}

\begin{rmk} The connection between the sets $Q_g$ and $P_g$ which we will present below is due to Brock and Granville and can be found in an early version of \cite{Brock}. There the connection is used to count the number of hyperelliptic curves in even characteristic, which is $a_0|_{g,even}$ in our terminology.\end{rmk}

\begin{lem} \label{lem-reform} Let $h$ and $f$ be polynomials over $k$. For any irreducible polynomial $m$ over $k$, the following two statements are equivalent:
\begin{itemize}
\item[$\star$] $m|\gcd(h,{f'}^2+f{h'}^2)$; 
\item[$\star$] there is a polynomial $l$ over $k$,  
such that $m|h$ and $m^2 | f+hl+l^2$. 
\end{itemize}
\end{lem}
\begin{pf}
Say that $\alpha \in k_n$ is a root of an irreducible polynomial $m$ and of the polynomial $\gcd(h,{f'}^2+f{h'}^2)$. Let $l$ be equal to $f^{q^n/2}$. Working modulo $(x-\alpha)^2$ we then get
\begin{multline*}f+hl+l^2 = f+hf^{q^n/2}+f^{q^n} \equiv f(\alpha)+ f'(\alpha)(x-\alpha)+h'(\alpha)f(\alpha)^{q^n/2}(x-\alpha)+f(\alpha)^{q^n} \\\equiv (x-\alpha)(f'(\alpha)+h'(\alpha)f(\alpha)^{q^n/2}) \equiv (x-\alpha)(f'(\alpha)^2+h'(\alpha)^2f(\alpha))^{1/2} = 0,\end{multline*}
which tells us that $m^2|f+hl+l^2$. For the other direction, assume that we have an irreducible polynomial $m$ and a polynomial $l$ such that $m|h$ and $m^2|f+hl+l^2$. Differentiating the polynomial $f+hl+l^2$ gives $m^2|f'+h'l+hl'$, and thus $m|f'+h'l$. Taking squares we get $m^2|f'^2+h'^2l^2$ and then it follows that $m^2|f'^2+h'^2(f+hl)$ and hence $m|f'^2+h'^2f$.\end{pf} 

Let $(h,f)$ be an element of $Q_g$. In the first part of the proof of Lemma~\ref{lem-reform}, we may take for $l$ any representative of $f^{q^n/2}$ modulo $h$, because for these $l$ we have $f+hl+l^2 \equiv  f+hf^{q^n/2}+f^{q^n}$ modulo $(x-\alpha)^2$. In the second part it does not matter which degree $l$ has. We conclude from this that Lemma~\ref{lem-reform} also holds if we assume that $l$ is of degree at most $g+1$.

Choose $g \geq -1$ and let $(h,f) \in Q_g$. Lemma~\ref{lem-reform} gives the following alternative formulation of the conditions \eqref{eq-deg}, \eqref{eq-nons} and \eqref{eq-inf}. For all polynomials $l$ of degree at most $g+1$:
\begin{gather}
\label{eq-nons2} m|h, \; m^2|f+hl+l^2 \implies \deg(m)=0; \\
\label{eq-inf2} \deg(h)=g+1 \quad \text{or} \quad \deg(f+hl+l^2) \geq 2g+1.
\end{gather}
Here we used that $t | \gcd(h_{\infty},f_{\infty}'^2+f_{\infty}h_{\infty}'^2)$ if and only if $t|h_{\infty}$ and there exists a polynomial $l_{\infty}$ such that $\deg(l_{\infty})\leq g+1$ and $t^2|f_{\infty}+h_{\infty}l_{\infty}+l_{\infty}^2$. In turn, this happens if and only if $\deg(h) \leq g$ and there exists a polynomial $l$ of degree at most $g+1$ such that $\deg(f+hl+l^2) \leq 2g$, where we connect $l$ and $l_{\infty}$ using the definitions $l:=x^{g+1}l_{\infty}(1/x)$ and $l_{\infty}:=t^{g+1}l(1/t)$.
 
This reformulation leads us to making the following definition.
\begin{dfn} 
Let $\sim_g$ be the relation on $Q_g$ given by $(h,f)\sim_g(h,f+hl+l^2)$ if $l$ is a polynomial of degree at most $g+1$. This is an equivalence relation and since $(h,f)=(h,f+hl+l^2)$ if and only if $l=0$ or $l=h$, the number of elements of each equivalence class $[(h,f)]_g$ is $q^{g+2}/2$. If $(h,f) \in P_g \subset Q_g$ then $[(h,f)]_g \subset P_g$ and we get an induced equivalence relation on $P_g$ which we also denote $\sim_g$. 
\end{dfn}

We will now construct all $\sim_g$ equivalence classes of elements of $Q_g$ in terms of the $\sim_i$ equivalence classes of the elements in $P_i$, where $i$ is between $-1$ and $g$. This is the counterpart of factoring a polynomial into a square-free part and a squared part in the case of odd characteristic.

\begin{dfn} For $z:=[(h,f)]_i \in P_i/\sim_i$ let $V_z$ be the set of all equivalence classes $[(mh,m^2f)]_{g}$ in $Q_g$ for all monic polynomials $m$ of degree at most $g-i$. This is well defined since if $(h_1,f_1) \sim_i (h_2,f_2)$ then $(mh_1,m^2f_1) \sim_g (mh_2,m^2f_2)$. 
\end{dfn}

\begin{lem} \label{lem-disjoint}
The sets $V_z$ for all $z \in P_i/\sim_i$ where $-1 \leq i \leq g$ are disjoint.
\end{lem}
\begin{pf}
Say that for some $z_1$ and $z_2$ the intersection $V_{z_1} \cap V_{z_2}$ is nonempty. That is, there exist $(h_1,f_1) \in P_{i_1}$, $(h_2,f_2) \in P_{i_2}$ and monic polynomials $m_1$, $m_2$ such that $m_1 h_1 = m_2 h_2$ and $m_1^2 f_1 = m_2^2f_2+m_2h_2l+l^2$. If for some irreducible polynomial $r$ we have $r|m_1$ but $r \nmid m_2$, it follows that $r|h_2$ and $r^2|m_2^2f_2+m_2h_2l+l^2$. By the equivalence of conditions \eqref{eq-nons} and \eqref{eq-nons2}, this implies that $r|(m_2^2f_2)'^2+m_2^2f_2(m_2h_2)'^2$ which in turn implies that $r|f_2'^2+f_2h_2'^2$. Since $(h_2,f_2) \in P_{i_2}$ we see that $r$ must be constant.  Hence every irreducible factor of $m_1$ is a factor of $m_2$. The situation is symmetric and therefore the converse also holds.

So far we have not ruled out the possibility that a factor in $m_1$ appears with higher multiplicity than in $m_2$, or vice versa. Let $m$ be the product of all irreducible factors of $m_1$ and put $\tilde m_1:=m_1/m$, $\tilde m_2:=m_2/m$ and $\tilde l :=l/m$. We are then in the same situation as above, that is $\tilde m_1 h_1 = \tilde m_2 h_2$ and $\tilde m_1^2 f_1 = \tilde m_2^2f_2+\tilde m_2 h_2 \tilde l+\tilde l^2$. Thus, if $r$ is an irreducible polynomial such that $r|\tilde m_1$ but $r \nmid \tilde m_2$ we can argue as above to conclude that $r$ is constant. By a repeated application of this line of reasoning we can conclude that $m_1$ and $m_2$ must be equal.

It now follows that $h_1=h_2$ and that $m_2|l$, thus $(h_1,f_1) \sim_{i_1} (h_2,f_2)$. This tells us that $V_{z_1} \cap V_{z_2}$ is only nonempty when $z_1=z_2$.\end{pf}

\begin{lem}  \label{lem-cover} 
The sets $V_z$ for all $z \in P_i/\sim_i$ where $-1 \leq i \leq g$ cover $Q_g/\sim_g$.
\end{lem}
\begin{pf} Pick any element $(h_1,f_1) \in Q_{g}$ and put $g_1:=g$. We define a procedure, where at the $i$th step we ask if there are any polynomials $m_i$ and $l_i$ such that $\deg(m_i) > 0$, $\deg(l_i) \leq g_i+1$, $m_i|h_i$ and $m_i^2|f_i+h_il_i+l_i^2$. If so, take any such polynomials $m_i$, $l_i$ and define $h_{i+1}:=h_i/m_i$, $f_{i+1}:=(f_i+h_il_i+l_i^2)/m_i^2$ and $g_{i+1}:=g_i-\deg(m_i)$. This procedure will certainly stop. Assume that the procedure has been carried out in some way and that it has stopped at the  $j$th step, leaving us with some pair of polynomials $(h_j,f_j)$.

Next, we take $(h_j,f_{j+1})$ to be any element of the set $[(h_j,f_j)]_{g_j}$ for which $\deg(f_{j+1})$ is minimal. Say that $f_{j+1}=f_j+h_jl_j+l_j^2$ where $\deg(l_j) \leq g_j+1$ and let us define $g_{j+1}$ to be the number such that $2g_{j+1}+1 \leq \max\bigl(2 \deg(h_j),\deg(f_{j+1})\bigr) \leq 2g_{j+1}+2$.
The claim is now that $(h_j,f_{j+1}) \in P_{g_{j+1}}$. By definition, condition \eqref{eq-deg} holds for $(h_j,f_{j+1})$. If there were polynomials $m_{j+1}$ and $l_{j+1}$ such that $m_{j+1}|h_j$ and $m_{j+1}^2|f_{j+1}+h_jl_{j+1}+l_{j+1}^2$ then the pair of polynomials $m_{j+1}$ and $l_j+l_{j+1}$ would contradict that the process above stopped at the $j$th step. Hence condition \eqref{eq-nons2} is fulfilled for $(h_j,f_{j+1})$. Condition \eqref{eq-inf2} is fulfilled if $2 \deg(h_j) \geq\deg(f_{j+1})$ because then $\deg(h_j)=g_{j+1}+1$. On the other hand, if $2 \deg(h_j) < \deg(f_{j+1})$ and there were a polynomial $l_{j+1}$ such that $\deg(l_{j+1}) \leq g_{j+1}+1$ and $\deg(f_{j+1}+h_jl_{j+1}+l_{j+1}^2) \leq 2g_{j+1}$ then this would contradict the minimality of $\deg(f_{j+1})$. We conclude that $(h_j,f_{j+1}) \in P_{g_{j+1}}$.

Finally we see that if we put $\hat m_r:=\prod_{i=1}^{r-1} m_{i}$ and $l:=\sum_{i=1}^{j} \hat m_i l_i$, then $\deg(l) \leq g+1$, $h_1=\hat m_j h_j$ and $f_1=\hat m_j^2 f_{j+1}+h_1l+l^2$. This shows that $V_z$ contains $[(h_1,f_1)]_g$ where $z:=[(h_j,f_{j+1})]_{g_{j+1}} \in P_{g_{j+1}}/\sim_{g_{j+1}}$.\end{pf}

Using the lemmas above we will be able to write $\UU_g$ in terms of $u_i$ for $i$ between $-1$ and $g$. After this we will determine $\UU_g$ for large enough values of $g$. We divide into two cases. 

\begin{ntn} Let $S_j$ denote all polynomials of degree at most $j$, and let $S'_j \subset S_j$ consist of the monic polynomials.\end{ntn}

\subsection{The case $\alpha \in A'(\nn)$} \label{sec-n2even} Fix an element $\alpha \in A'(\nn)$. It follows from Lemma~\ref{lem-disjoint} and Lemma~\ref{lem-cover} that
\begin{equation}  \label{eq-Uua-even} \UU_{g,\alpha} = \frac{I_g}{2}  
\sum_{l \in S_{g+1}} \sum_{j=-1}^{g}  \sum_{z \in P_{j}/\sim_{j}} \sum_{[(h,f)]_j \\ \in V_z} \prod_{i=1}^{m} \tau_{n_i} \bigl(h(\alpha_{i}), (f+h l+l^2)(\alpha_{i}) \bigr)^{r_i}. \end{equation}

\begin{lem} \label{lem-tauprop} Choose any $s \geq 1$ and $t_1,t_2$ in $k_s$. We then have
\begin{gather*}
\label{eq-tau1} \tau_s(vt_1,v^2t_2)=\tau_s(t_1,t_2) \quad \text{for all $v\neq 0 \in k_s$;} \\
\label{eq-tau2} \tau_s(t_1,t_2+vt_1+v^2)=\tau_s(t_1,t_2) \quad \text{for all $v \in k_s$}. 
\end{gather*} 
\end{lem}
\begin{pf} Clear. \end{pf}

Fix elements $z=[(h_0,f_0)]_i \in P_i/\sim_i$ and $\beta \in \Ab^1(k_s)$ and define $V'_z$ to be the subset of $V_z$ of classes $[(\tilde m h_0,\tilde m^2f_0)]_g$, where $\tilde m$ is a monic polynomial with $\tilde m(\beta) \neq 0$. Lemma~\ref{lem-tauprop} shows that $\tau_s(h(\beta),f(\beta))$ is constant for all $s$ and $(h,f)$ such that $[(h,f)]_g \in V'_z$. Applying this to equation \eqref{eq-Uua-even} after recalling Definition \ref{dfn-b} we find that
\begin{multline}  \label{eq-receven} \UU_{g,\alpha} = I_g\frac{ q^{g+1}}{2} \sum_{j=-1}^{g}  \sum_{z \in P_{j}/\sim_{j}} \sum_{\tilde m \in S'_{g-j}}  \prod_{i=1}^{m} \tau_{n_i} \bigl((\tilde m h)(\alpha_{i}), (\tilde m^2f)(\alpha_{i}) \bigr)^{r_i} =\\= I_g  \frac{ q^{g+1}}{2} \sum_{j=-1}^{g}\bb_{g-j} u_{j,\alpha} \frac{2}{q^{j+1}I_{j}}=\sum_{i=0}^{g+1}\bb_{i} u_{g-i,\alpha},  \end{multline} 
where we have taken into account that the group of isomorphisms depends upon $g$ and that the numbers of elements of the equivalence classes of the relations $\sim_{g-j}$ and $\sim_g$ differ by a factor $q^{j}$. From the definitions we see that $q^{g-j} \, I_{g}/I_{j}=1$.

For any $g \geq -1$ and any $h_0 \in S_{g+1}$ it is clear that
\begin{equation}\label{eq-UUr2even} \sum_{(h_0,f) \in Q_g}\prod_{i=1}^m \tau_{n_i} \bigl(h_0(\alpha_i),f(\alpha_{i}) \bigr)^{r_{i}} = \begin{cases}0 & \text{if $\forall i:r_i=2$, $\exists j:h_0(\alpha_j)=0$};  \\q^{2g+3} & \text{if $\forall i:r_i=2$, $\forall j: h_0(\alpha_j) \neq 0$}. \end{cases}\end{equation}

For any $g$ such that $2g+2 \geq \abs{\nn}-1$, and any nonzero polynomial $h_0$ of degree at most $g+1$, Lemma~\ref{lem-onetoone} tells us that
\begin{multline}\label{eq-UUr1even} \sum_{(h_0,f_1+p_{\alpha} f_2) \in Q_g}\prod_{i=1}^m \tau_{n_i} \bigl(h_0(\alpha_i),(f_1+p_{\alpha}f_2 )(\alpha_{i}) \bigr)^{r_{i}} = \\ = q^{2g+3-\abs{\nn}} \sum_{f_1 \in S_{\abs{\nn}-1}} \prod_{i=1}^{m} \tau_{n_i} \bigl( h_0(\alpha_{i}), f_1(\alpha_{i}) \bigr)^{r_i} = \\ = q^{2g+3-\abs{\nn}} \sum_{(\beta_1,\ldots,\beta_m) \in \prod_{i=1}^m k_{n_i} } \prod_{i=1}^{m} \tau_{n_i} \bigl( h_0(\alpha_{i}), \beta_i \bigr)^{r_i} = 0 \quad \text{if $\exists i:r_i=1$}, 
\end{multline}
because for all $a\in k_s$ there are as many $b \in k_s$ for which $\tau_s(a,b)=1$ as there are $b \in k_s$ for which $\tau_s(a,b)=-1$. 

Summing equations \eqref{eq-UUr2even} and \eqref{eq-UUr1even} over all $h_0 \in S_{g+1}$ and using that $q^{2g+3}I_g=I q^{g+1}$ we get
\begin{equation} \label{eq-UUeven}
\UU_{g,\alpha} = \begin{cases} I\, (q-1) q^{g+1} \bb_{g+1}  & \text{if $\forall i:r_i=2$, $g \geq -1$;} \\0 & \text{if $\exists i:r_i=1$, $g \geq \frac{\abs{\nn}-3}{2}$.} \end{cases}
\end{equation}

\subsection{The case $\alpha \in A(\nn) \setminus A'(\nn)$} \label{sec-n1even} Fix an $\alpha \in A(\nn)\setminus A'(\nn)$. We can assume that $\alpha_1 = \infty$, and then $\tilde \alpha:=(\alpha_2,\ldots,\alpha_m) \in A'(\tilde \nn)$ where $\tilde \nn:=(n_2,\ldots,n_m)$. 

\begin{lem} \label{lem-tauinft}
For any element $(h,f)\in P'_i$ and any monic polynomial $m$ of degree $g-i$,
\begin{gather*}\tau_s((mh)(\infty),(m^2f)(\infty))=\tau_s(h(\infty),f(\infty)); \\ \tau_s((mh)(\infty),(f+lh+l^2)(\infty))=\tau_s(h(\infty),f(\infty)).\end{gather*} 
\end{lem}
\begin{pf} Clear. \end{pf}

For any $(h,f) \in Q_g$ it holds that if $\mathrm{deg}(h) < g+1$ then $\tau_s(h(\infty),f(\infty))=0$ for all $s$. Define therefore $P'_g$ and $Q'_g$ to be the subsets of $P_g$ and $Q_g$ respectively, that consist of pairs $(h,f)$ such that $\mathrm{deg}(h)=g+1$. We get an induced relation $\sim_i$ on $P'_i$ and $Q'_i$ and we let $V''_z$ be the set of all equivalence classes $[(mh,m^2f)]_{g}$ in $Q'_g$ for all monic polynomials $m$ of degree $g-i$, where $z:=[(h,f)]_i \in P'_i/\sim_i$. In the same way as in Lemma~\ref{lem-disjoint} and \ref{lem-cover} we see that the sets $V''_z$ for all $z \in P'_i/\sim_i$, where $-1 \leq i \leq g$, are disjoint and cover $Q'_g/\sim_g$. 
Using this together with Lemma~\ref{lem-tauinft} and the arguments showing equation~\eqref{eq-receven} we find that
\begin{multline}  \label{eq-receveninf} \UU_{g,\alpha} = \frac{I_g}{2}  \sum_{l \in S_{g+1}} \sum_{z \in Q'_{g}/\sim_{g}} \prod_{i=1}^{m} \tau_{n_i} \bigl(h(\alpha_{i}), (f+hl+l^2)(\alpha_{i}) \bigr)^{r_i} = \\ 
= I_g  \frac{q^{g+1}}{2} \sum_{j=-1}^{g}  \sum_{z \in P'_{j}/\sim_{j}} \sum_{\tilde m \in R'_{g-j}}  \prod_{i=1}^{m} \tau_{n_i} \bigl((\tilde m h)(\alpha_{i}), (\tilde m^2 f)(\alpha_{i}) \bigr)^{r_i} =\sum_{i=0}^{g+1} b^{\tilde \nn}_{i} u_{g-i,\alpha}.
\end{multline} 

If we choose $g$ such that $2g+2 \geq \abs{\nn} -1$, $h_0 \in R_{g+1}$ and we put $p_{\alpha}(x):=x \, p_{\tilde \alpha}$, then we find in the same way as for equation~\eqref{eq-UUr1even} that
\begin{multline} \label{eq-UUr1infeven} \sum_{(h_0,f_1+p_{\alpha} f_2) \in Q_g}\prod_{i=1}^m \tau_{n_i} \bigl(h_0(\alpha_i),(f_1+p_{\alpha}f_2 )(\alpha_{i}) \bigr)^{r_{i}} =\\ = q^{2g+3-\abs{\nn}} \sum_{(\beta_1,\ldots,\beta_m) \in \prod_{i=1}^m k_{n_i} } \prod_{i=1}^{m} \tau_{n_i} \bigl( h_0(\alpha_{i}), \beta_i \bigr)^{r_i} = 0 \quad \text{if $\exists i:r_i=1$}.
\end{multline}
Since equation~\eqref{eq-UUr2even} also hold for $\alpha \in A(\nn) \setminus A'(\nn)$ we find, by summing over all polynomials $h_0 \in R_{g+1}$, that
\begin{equation} \label{eq-Ueven}
\UU_{g,\alpha} = \begin{cases} I \, (q-1)  q^{g+1} b_{g+1}^{\tilde \nn}  & \text{if $\forall i:r_i=2$, $g \geq -1$;} \\0 & \text{if $\exists i:r_i=1$, $g \geq \frac{\abs{\nn}-3}{2}$.} \end{cases}
\end{equation}

\subsection{The two cases joined}
Recall that $J=(q-1) \,I \,\abs{A(\nn)}$.

\begin{thm} \label{thm-rec1even} For any pair $(\nn;\rr) \in \mathcal{N}_m$,
\begin{equation*} \sum_{j=0}^{g+1}\bb_{j} u_{g-j} =\begin{cases} J \, q^{g+1} \bb_{g+1}  & \text{if $\forall i:r_i=2$, $g \geq -1$}; \\0 & \text{if $\exists i:r_i=1$, $g \geq \frac{\abs{\nn}-3}{2} $}.\end{cases}
\end{equation*} 
\end{thm}
\begin{pf} The theorem follows from combining equations \eqref{eq-receven}, \eqref{eq-UUeven}, \eqref{eq-receveninf} and \eqref{eq-Ueven}, using Lemma~\ref{lem-b}.\end{pf}

\begin{thm} \label{thm-rec2even} For any pair $(\nn;\rr) \in \mathcal{N}_m$,
\begin{equation*} \sum_{j=0}^{\min(\abs{\nn}-1,g+1)}(\bb_{j}-q\bb_{j-1}) u_{g-j} = \begin{cases} J \, q^{g+1} (\bb_{g+1}-\bb_{g})  & \text{if $\forall i:r_i=2$, $g \geq 0$}; \\0 & \text{if $\exists i:r_i=1$, $g \geq \frac{\abs{\nn}-1}{2}$}. \end{cases}
\end{equation*} 
\end{thm}
\begin{pf} In the notation of the proof of Theorem~\ref{thm-rec2}, the theorem follows from applying Theorem~\ref{thm-rec1even} to the expression $F(g)-q F(g-1)$.\end{pf}

\begin{thm} \label{thm-chareqeven} By applying Theorem~\ref{thm-rec2even} to each pair $(\nn;\rr)$ appearing in the decomposition (given by Lemma~\ref{lem-decompeven}) of $a_{\lambda}|_{g,even}$ we get a linear recurrence relation for $a_{\lambda}|_{g,even}$. The characteristic polynomial of this linear recurrence relation equals~\eqref{eq-chareq}.
\end{thm}
\begin{pf} We know that the decomposition of $a_{\lambda}|_g$ is independent of characteristic, and since the left hand side of the equation in Theorem~\ref{thm-rec2even} is the same as the left hand side of the equation of Theorem~\ref{thm-rec2} this theorem follows in the same way as Theorem~\ref{thm-chareq}. \end{pf}

\section{Results for weight up to $7$ in even characteristic}\label{sec-results2}
In this section we compute, for any number $g$ and any finite field $k$ of even characteristic, all $a_{\lambda}|_{g,even}$ of weight at most $7$. First we will exploit the similarities of Theorems~\ref{thm-rec1} and~\ref{thm-rec1even}.

\begin{lem} \label{lem-tocomp} If $g \geq n-2$ then $\bb_{2g+2}=q^{g+1}\bb_{g+1}$. \end{lem}
\begin{pf} Fix a pair $(\nn;\rr) \in \mathcal{N}_m$. Lemma~\ref{lem-bcoeff} tells us that $\bb_j=q\bb_{j-1}+d_{\abs{\nn}-1-j}$, so if $j \geq \abs{\nn}$ then $\bb_j=q\bb_{j-1}$ and thus $\bb_j=q^{j+1-\abs{\nn}}\bb_{\abs{\nn}-1}$.\end{pf}

\begin{rmk} \label{rmk-compare} If $r_i=1$ for some $i$ and $g \geq (\abs{\nn}-3)/2$, then the recursive relations of Theorems~\ref{thm-rec1even} and~\ref{thm-rec1} are equal. On the other hand, if $r_i=2$ for all $i$ we see from Lemma~\ref{lem-tocomp} that the recursive relations of Theorems~\ref{thm-rec1even} and~\ref{thm-rec1} are equal if $g \geq \abs{\nn}-2$.
\end{rmk}

\begin{thm} \label{thm-ind} For weight less than or equal to $5$, $a_{\lambda}|_{g,even} = a_{\lambda}|_{g,odd}$ as functions (in this case polynomials) in $q$.
\end{thm}
\begin{pf} Consider any $a_{\lambda}|_g$ with $\abs{\lambda} \leq 5$. By Lemma \ref{lem-decr} it suffices to show that $u_g$ is independent of characteristic when $(\nn;\rr)\in \mathcal{N}_m$ is such that $\sum_{i=1}^m n_ir_i \leq 5$. Clearly $u_{-1}=J$ is always independent of characteristic. Clearly, Lemma~\ref{lem-g0} also holds in even characteristic. We can therefore assume that $r_i=2$ for all $i$ in the case of genus $0$. But if $r_i=2$ for all $i$ then $\abs{\nn} \leq 2$ and hence, by Remark~\ref{rmk-compare}, $u_0$ will be independent of characteristic. 

This takes care of the base cases of the recurrence relations for $u_g$ when $g \geq 1$, given by Theorems \ref{thm-rec1} and \ref{thm-rec1even}. Again by Remark~\ref{rmk-compare} we see that (both in the case when $r_i=2$ for all $i$, and when $r_i=1$ for some $i$) when $g \geq 1$ these recurrence relations are the same. We can therefore conclude that $u_g$ is independent of characteristic for all $g$. \end{pf} 

We will now compute $a_{\lambda}|_{g,even}$ for weight $6$ in the same way as in Section~\ref{sec-wsix}. To compute $u_g$ of degree at most $5$ using Theorem \ref{thm-rec1even} we need to find the base case $u_0$.  But when the genus is $0$ we can use Lemma~\ref{lem-g0} (which also holds in even characteristic) to reduce to the case that $r_i=2$ for all $i$, which is always computable using Theorem \ref{thm-rec1even}.

What is left is the general case of the decomposition of $a_{\lambda}|_{g,even}$. We then need the base cases of genus $0$ and $1$. Again, the genus $0$ part is no problem. The computation of $a_{\lambda}|_1$ in \cite{Jonas2} is independent of characteristic. We can therefore compute the genus $1$ part (compare Section~\ref{sec-wsix}).

\begin{rmk}  As in the case of odd characteristic, for all $g$ and all $\lambda$ such that $\abs{\lambda} \leq 7$, $a_{\lambda}|_{g,even}$ is a polynomial when considered as a function in the number $q$ (compare Remark~\ref{rmk-polynomial}) of elements of the finite field $k$ of even characteristic.

In Theorem~\ref{thm-ind} we saw that the polynomial functions $a_{\lambda}|_{g,odd}$ and $a_{\lambda}|_{g,even}$ are equal (for a fixed $g$), if $\abs{\lambda} \leq 5$. But for weight $6$ there are $\lambda$ such that the two polynomials are different, this occurs for the first time for genus $3$, see Example~\ref{exa-a16even}. 
\end{rmk}

\begin{exa} We wish to compute $u_{g,even}$ in the case $(\nn;\rr)=((1,1);(2,2,2))$. We see that $u_{-1}=1$ and Theorem~\ref{thm-rec1even} gives $u_0=q^2-3q+2$. This result is different from the $1$ in the case of odd characteristic, see Example \ref{exa-u2}. Continued use of Theorem~\ref{thm-rec1even} gives $u_1=q^4-3q^3+5q^2-6q+3$ and then Theorem~\ref{thm-rec2even} gives
$$u_g=2u_{g-1} - u_{g-2} + q^{2g-1}(q-1)^3 \quad \text{for $g \geq 2$}.$$
Solving this leaves us with
$$u_{g,even}^{((1,1,1);(2,2,2))}=\frac{(q-1)(q^{2g+3}+g(q^2-1)-3q-2)}{(q+1)^2}.$$
\end{exa}

\begin{exa} \label{exa-a16even} The result for $a_{[1^6]}|_{g,even}$ is
\begin{equation*}
a_{[1^6]}|_{g,even} = a_{[1^6]}|_{g,odd}-\frac{5}{8} \, g(g-1)(g-2)\bigl((g-3)(q-1)-4\bigr).
\end{equation*}
\end{exa}

\begin{exa} The result for $a_{[1^2,4]}|_{g,even}$ is
\begin{equation*}
a_{[1^2,4]}|_{g,even}=a_{[1^2,4]}|_{g,odd}-\frac{1}{4}\begin{cases}g(q-1) & \text{if $g\equiv 0$ mod $4$};\\(g-1)(q-1) & \text{if $g\equiv 1$ mod $4$};\\(g-2)(q-1) & \text{if $g\equiv 2$ mod $4$};\\(g-3)(q-1)-4 & \text{if $g\equiv 3$ mod $4$}. \end{cases}
\end{equation*}
\end{exa}

\section{Cohomological results} \label{sec-coh} 
\subsection{Cohomological results for $\Hh_{g,n}$} \label{sec-coh-et}
Define the local system $\V:=R^1\pi_{*}(\Ql)$ where $\pi : \mathcal{M}_{g,1} \to \mathcal{M}_{g}$ is the universal curve. For every partition (note that in this section we use a different notation for partitions) $\lambda=(\lambda_1 \geq \ldots \geq \lambda_g \geq 0)$ there is an irreducible representation of $\mathrm{GSp}(2g)$ with highest weight $(\lambda_1-\lambda_2)\gamma_1+\ldots+\lambda_g \gamma_g-|\lambda| \eta$, where the $\gamma_i$ are suitable fundamental roots and $\eta$ is the multiplier representation, and we define $\Vla$ to be the corresponding local system. Let us also denote by $\Vla$ its restriction to $\Hh_g$. In Lemma \ref{lem-loca} below we will see that making an $\s_{\tilde n}$-equivariant count of points of $\Hh_{g,\tilde n}$ over a finite field $k$, for all $\tilde n \leq n$, is equivalent to computing the trace of Frobenius on the compactly supported $\ell$-adic Euler characteristic $\Eulc(\Hh_{g} \otimes \bar k,\Vla)$, for every $\lambda$ with $\abs{\lambda} \leq n$ (where $\ell\nmid \abs{k}$). For more details, see \cite{G-2} and \cite{G-res}.

Thus, we can use the results of Section \ref{sec-results} together with Theorem 3.2 in \cite{Jonas2} to compute the $\ell$-adic Euler characteristic $\Eulc(\Hh_{g} \otimes \overline{\Q},\Vla)$ in $K_0(\mathsf{Gal}_{\Q})$, the Grothen\-dieck group of $\mathrm{Gal(\bar{\Q}/\Q})$-representations, for every $\lambda$ with $\abs{\lambda} \leq 7$. Specifically, Theorem~3.2 in \cite{Jonas2} tells us that if there is a polynomial $P$ such that $\mathrm{Tr}(F,\Eulc(\Hh_{g} \otimes \bar k,\Vla))=P(q)$ for all finite fields $k$, possibly with the exception of a finite number of characteristics, then $\Eulc(\Hh_{g} \otimes \overline{\Q},\Vla)=P(\mathbf{q})$, where $\mathbf{q}$ is the class of $\Ql(-1)$ in $K_0(\mathsf{Gal}_{\Q})$. By excluding even characteristic, Section \ref{sec-results} (see Remark~\ref{rmk-polynomial}) and Lemma \ref{lem-loca} shows that there is indeed such a polynomial for all $g$ and all $\abs{\lambda} \leq 7$.

\begin{exa} For $g=8$ and $\lambda=(5,1)$ we have
$$\Eulc(\Hh_{g} \otimes \overline{\Q},\Vla)=5\mathbf{q}^5-28\mathbf{q}^4+4\mathbf{q}^3+96\mathbf{q}^2-34\mathbf{q}-\mathbf{88}.$$
\end{exa}

\subsection{Cohomological results for $\Mmb_{2,n}$ and $\Mm_{2,n}$}  \label{sec-coh-bet}
Using the stratification of $\Mmb_{g,n}$ we can make an $\s_n$-equivariant count of its number of points using the $\s_n$-equivariant counts of the points of $\Mm_{\tilde{g},\tilde{n}}$ for all $\tilde{g} \leq g$ and $\tilde{n}\leq n+2(g-\tilde{g})$ (see \cite[Thm~8.13]{GK} and also \cite{Mbar4}). Since all curves of genus $2$ are hyperelliptic, $\Mm_{2,n}$ is equal to $\Hh_{2,n}$. Above, we have made $\s_n$-equivariant counts of $\Hh_{2,n}$ for $n \leq 7$ and they were all found to be polynomial in $q$. These $\s_n$-equivariant counts can now be complemented with ones of $\Mm_{1,n}$ for $n \leq 9$ (see \cite[Section 15]{Jonas2}) and of $\Mm_{0,n}$ for $n \leq 11$ (see \cite[Prop 2.7]{Lehrer}), which are also found to be polynomial in $q$. We can then apply Theorem 3.4 in \cite{Mbar4} to conclude, for all $n\leq 7$, the $\s_n$-equivariant $\mathsf{Gal}_{\Q}$ (resp. Hodge) structure of the $\ell$-adic (resp. Betti) cohomology of $\Mmb_{2,n}$.

In the theorems below we give the $\s_n$-equivariant Hodge Euler characteristic (which by purity is sufficient to conclude the Hodge structure) in terms of the Schur polynomials and $\Ll$, the class of the Tate Hodge structure of weight $2$ in $K_0(\mathsf{HS}_\Q)$, the Grothen\-dieck group of rational Hodge structures. That is,
the action of $\s_{n}$ on $\Mmb_{2,n}$ induces an action on its cohomology, and hence $H^i(\Mmb_{2,n} \otimes \C,\Q)$ may be written as a direct sum of $H^i_{\lambda}(\Mmb_{2,n}\otimes \C,\Q)$, which correspond to the irreducible representations of $\s_n$ indexed by $\lambda \vdash n$ and with characters $\chi_{\lambda}$. In terms of this, 
the coefficient of the Schur polynomial $s_{\lambda}$ is equal to $1/\chi_{\lambda}(id) \cdot \sum_i (-1)^i[H^i_{\lambda}(\Mmb_{2,n}\otimes \C,\Q)] $. The results for $n \leq 3$ were previously known by the work of Getzler in \cite[Section 8]{G-2}. 

\begin{thm}The $\s_n$-equivariant Hodge Euler characteristic of $\Mmb_{2,4}$ is equal to
\begin{gather*}
(\Ll^7+8\Ll^6+33\Ll^5+67\Ll^4+67\Ll^3+33\Ll^2+8\Ll+\mathbf{1})s_{4}\\
+(4\Ll^6+26\Ll^5+60\Ll^4+60\Ll^3+26\Ll^2+4\Ll)s_{31}\\
+(2\Ll^6+12\Ll^5+28\Ll^4+28\Ll^3+12\Ll^2+2\Ll)s_{2^2}\\
+(3\Ll^5+10\Ll^4+10\Ll^3+3\Ll^2)s_{21^2}
\end{gather*}
\end{thm}

\begin{thm}The $\s_n$-equivariant Hodge Euler characteristic of $\Mmb_{2,5}$ is equal to
\begin{gather*}
(\Ll^8+9\Ll^7+49\Ll^6+128\Ll^5+181\Ll^4+128\Ll^3+49\Ll^2+9\Ll+\mathbf{1})s_{5}\\
+(6\Ll^7+48\Ll^6+156\Ll^5+227\Ll^4+156\Ll^3+48\Ll^2+6\Ll)s_{41}\\
+(3\Ll^7+31\Ll^6+106\Ll^5+159\Ll^4+106\Ll^3+31\Ll^2+3\Ll)s_{32}\\
+(8\Ll^6+42\Ll^5+65\Ll^4+42\Ll^3+8\Ll^2)s_{31^2}\\
+(6\Ll^6+26\Ll^5+43\Ll^4+26\Ll^3+6\Ll^2)s_{2^21}\\
+(\Ll^5+3\Ll^4+\Ll^3)s_{21^3}\end{gather*}
\end{thm}

\begin{thm}The $\s_n$-equivariant Hodge Euler characteristic of $\Mmb_{2,6}$ is equal to
\begin{gather*}
(\Ll^9+11\Ll^8+68\Ll^7+229\Ll^6+420\Ll^5+420\Ll^4+229\Ll^3+68\Ll^2+11\Ll+\mathbf{1})s_{6}\\
+(7\Ll^8+75\Ll^7+317\Ll^6+641\Ll^5+641\Ll^4+317\Ll^3+75\Ll^2+7\Ll)s_{5 1}\\
+(5\Ll^8+62\Ll^7+292\Ll^6+615\Ll^5+615\Ll^4+292\Ll^3+62\Ll^2+5\Ll)s_{4 2}\\
+(\Ll^8+21\Ll^7+108\Ll^6+236\Ll^5+236\Ll^4+108\Ll^3+21\Ll^2+\Ll)s_{3^2}\\
+(17\Ll^7+118\Ll^6+278\Ll^5+278\Ll^4+118\Ll^3+17\Ll^2)s_{4 1^2}\\
+(16\Ll^7+115\Ll^6+277\Ll^5+277\Ll^4+115\Ll^3+16\Ll^2)s_{3 2 1}\\
+(3\Ll^7+22\Ll^6+53\Ll^5+53\Ll^4+22\Ll^3+3\Ll^2)s_{2^3}\\
+(9\Ll^6+29\Ll^5+29\Ll^4+9\Ll^3)s_{31^3}\\
+(6\Ll^6+21\Ll^5+21\Ll^4+6\Ll^3)s_{2^2 1^2}
\end{gather*}
\end{thm}

\begin{thm}The $\s_n$-equivariant Hodge Euler characteristic of $\Mmb_{2,7}$ is equal to
\begin{gather*}
(\Ll^{10}+12\Ll^9+90\Ll^8+363\Ll^7+854\Ll^6+1125\Ll^5+854\Ll^4+363\Ll^3+90\Ll^2+\ldots)s_{7}\\ 
+(9\Ll^9+109\Ll^8+580\Ll^7+1529\Ll^6+2109\Ll^5+1529\Ll^4+580\Ll^3+109\Ll^2+9\Ll)s_{6 1}\\
+(6\Ll^9+100\Ll^8+606\Ll^7+1728\Ll^6+2430\Ll^5+1728\Ll^4+606\Ll^3+100\Ll^2+6\Ll)s_{5 2}\\
+(3\Ll^9+58\Ll^8+389\Ll^7+1153\Ll^6+1647\Ll^5+1153\Ll^4+389\Ll^3+58\Ll^2+3\Ll)s_{4 3}\\
+(28\Ll^8+258\Ll^7+831\Ll^6+1221\Ll^5+831\Ll^4+258\Ll^3+28\Ll^2)s_{5 1^2}\\
+(34\Ll^8+331\Ll^7+1133\Ll^6+1675\Ll^5+1133\Ll^4+331\Ll^3+34\Ll^2)s_{4 2 1}\\
+(12\Ll^8+140\Ll^7+489\Ll^6+738\Ll^5+489\Ll^4+140\Ll^3+12\Ll^2)s_{3^2 1}\\
+(8\Ll^8+91\Ll^7+335\Ll^6+502\Ll^5+335\Ll^4+91\Ll^3+8\Ll^2)s_{3 2^2}\\
+(28\Ll^7+143\Ll^6+228\Ll^5+143\Ll^4+28\Ll^3)s_{4 1^3}\\
+(34\Ll^7+170\Ll^6+275\Ll^5+170\Ll^4+34\Ll^3)s_{3 2 1^2}\\
+(10\Ll^7+47\Ll^6+77\Ll^5+47\Ll^4+10\Ll^3)s_{2^3 1}\\
+(4\Ll^6+7\Ll^5+4\Ll^4)s_{3 1^4}\\
+(2\Ll^6+6\Ll^5+2\Ll^4)s_{2^2 1^3}
\end{gather*}
\end{thm}

In Table~\ref{tab-coh} we present the nonequivariant information (remember that all cohomology is Tate) in the form of Betti numbers of $\Mmb_{2,n}$ for all $n \leq 7$. Notice that the table only contains as many numbers as we need to be able to fill in the missing ones using Poincar\'e duality. These results agree with Table 2 of ordinary Euler characteristics for $\Mmb_{2,n}$ for $n \leq 6$ found in \cite{Bini-Harer}. 

\begin{table}[htbp] \caption{Dimensions of $H^i(\Mmb_{2,n}\otimes \C, \Q)$ for $n \leq 7$.} \label{tab-coh}
\centerline{
\vbox{
\offinterlineskip
\hrule
\halign{&\vrule#& \quod \hfil#\hfil \strut \quod \cr
height2pt&\omit&&\omit&&\omit&&\omit&&\omit&&\omit&&\omit& \cr 
&  && $H^0$ && $H^2$ && $H^4$ && $H^6$ && $H^8$ && $H^{10}$ & \cr
height2pt&\omit&&\omit&&\omit&&\omit&&\omit&&\omit&&\omit& \cr 
\noalign{\hrule}
height2pt&\omit&&\omit&&\omit&&\omit&&\omit&&\omit&&\omit& \cr 
&$\Mmb_{2}$  && 1 && 2  &&   &&   &&   &&  &\cr
height2pt&\omit&&\omit&&\omit&&\omit&&\omit&&\omit&&\omit& \cr 
&$\Mmb_{2,1}$  && 1 && 3  && 5  &&   &&   &&  &\cr
height2pt&\omit&&\omit&&\omit&&\omit&&\omit&&\omit&&\omit& \cr 
&$\Mmb_{2,2}$  && 1 && 6  && 14  &&   &&   &&  &\cr
height2pt&\omit&&\omit&&\omit&&\omit&&\omit&&\omit&&\omit& \cr 
&$\Mmb_{2,3}$  && 1 && 12  && 44  && 67  &&   &&  &\cr
height2pt&\omit&&\omit&&\omit&&\omit&&\omit&&\omit&&\omit& \cr 
&$\Mmb_{2,4}$  && 1 && 24  && 144  && 333  &&   &&  &\cr
height2pt&\omit&&\omit&&\omit&&\omit&&\omit&&\omit&&\omit& \cr 
&$\Mmb_{2,5}$  && 1 && 48  && 474  && 1668  && 2501  &&  &\cr
height2pt&\omit&&\omit&&\omit&&\omit&&\omit&&\omit&&\omit& \cr 
&$\Mmb_{2,6}$  && 1 && 96  && 1547  && 8256  && 18296  &&  &\cr
height2pt&\omit&&\omit&&\omit&&\omit&&\omit&&\omit&&\omit& \cr 
&$\Mmb_{2,7}$  && 1 && 192  && 4986  && 39969  && 129342 && 189289  &\cr
height2pt&\omit&&\omit&&\omit&&\omit&&\omit&&\omit&&\omit& \cr 
} \hrule}}
\end{table}

The theorem used above also gives the corresponding results for $\Mm_{2,n}$ for $n \leq 7$, which we will present in terms of local systems $\Vla$ defined as above, but starting from $\V:=R^1\pi_* \Q$. See~\cite[Section 8]{G-2} for the results on $\Eulc(\Mm_{2}\otimes \C,\Vla)$, for all $\lambda$ of weight at most $3$.

\begin{thm} \label{thm-M2}
The Hodge Euler characteristics of the local systems $\Vla$ on $\Mm_{2}':=\Mm_{2}\otimes \C $ of weight $4$ or $6$ are equal to
\begin{multline*}
\Eulc(\Mm_{2}',\V_{(4,0)})=\mathbf{0}, \;\; \Eulc(\Mm_{2}',\V_{(3,1)})=\Ll^2-\mathbf{1}, \;\; \Eulc(\Mm_{2}',\V_{(2,2)})=-\Ll^4, \\
\Eulc(\Mm_{2}',\V_{(6,0)})=-\mathbf{1},\;\;  \Eulc(\Mm_{2}',\V_{(5,1)})=\Ll^2-\Ll-\mathbf{1}, \\ \Eulc(\Mm_{2}',\V_{(4,2)})=\Ll^3,\;\;  \Eulc(\Mm_{2}',\V_{(3,3)})=-\Ll-\mathbf{1}.
\end{multline*}
\end{thm}

\section{Appendix: Introducing $b_i$, $c_i$ and $r_i$} \label{sec-app1}
This section will give an interpretation of the information carried by the $u_g$'s. It will be in terms of counts of hyperelliptic curves together with prescribed inverse images of points on $\Pb^1$ under their unique degree $2$ morphism. 

\begin{dfn}
Let $C_{\varphi}$ be a curve defined over $k$ together with a separable degree~$2$ morphism $\varphi$ over $k$ from $C$ to $\Pb^1$. We then define
\begin{gather*}
b_{i}(C_{\varphi}):=\abs{\{ \alpha \in A(i) : \abs{\varphi^{-1}(\alpha)}=2, \varphi^{-1}(\alpha) \subseteq C(k_i)  \}}, \\
c_{i}(C_{\varphi}):=\abs{\{ \alpha \in A(i) : \abs{\varphi^{-1}(\alpha)}=2, \varphi^{-1}(\alpha) \nsubseteq C(k_{i}) \}}
\end{gather*}
and put $r_i(C_{\varphi}):=b_{i}(C_{\varphi})+c_{i}(C_{\varphi})$.
\end{dfn}

The number of ramification points of $f$ that lie in $A(i)$ is then equal to $\abs{A(i)}-r_{i}(C_{\varphi}).$ Let $\lambda_i$ denote the partition of $i$ consisting of one element. We then find that 
$$\abs{C_{\varphi}(\lambda_i)}=\abs{A(i)}+b_{i}(C_{\varphi})-c_{i}(C_{\varphi})+\begin{cases} 2c_{i/2}(C_{\varphi}) & \text{if $i$ is even}; \\ 0 & \text{if $i$ is odd}. \end{cases} $$ 
and thus
$$a_n(C_{\varphi})=\sum_{i | n \,:\, 2i \nmid n} \bigl(c_i(C_{\varphi})-b_i(C_{\varphi})\bigr) + \sum_{i:2i | n} \bigl(-b_i(C_{\varphi})-c_i(C_{\varphi})\bigr).$$

\begin{dfn} \label{dfn-bc} For partitions $\mu$ and $\nu$, $g \geq 2$ and odd characteristic, define
\begin{equation*}  b_{\mu} c_{\nu}|_g:= \sum_{[C_f] \in \Hh_g(k)/\cong_k} \frac{1}{\abs{\Aut_k(C_f)}}\,\prod_{i=1}^{l(\mu)} b_{i}(C_f)^{\mu_i} \, \prod_{j=1}^{l(\nu)} c_{j}(C_f)^{\nu_j}.\end{equation*}
The number $\abs{\mu}+\abs{\nu}$ will be called the weight of this expression.
\end{dfn}

\begin{rmk} We can, in the obvious way, also define $a_{\lambda} b_{\mu} c_{\nu}|_g,$ but from the relation between $a_i(C_f)$, $b_i(C_f)$ and $c_i(C_f)$ we see that this gives no new phenomena.
\end{rmk}

Directly from the definitions we get the following lemma.
\begin{lem} Let the characteristic be odd and let $f$ be an element of $P_g$. We then have
$$b_{i}(C_f)=\frac{1}{2} \sum_{\alpha \in A(i)} \Bigl(\chi_{2,i}\bigl(f(\alpha)\bigr)^2+\chi_{2,i}\bigl(f(\alpha)\bigr)\Bigr)$$
and
$$c_{i}(C_f)=\frac{1}{2} \sum_{\alpha \in A(i)}\Bigl(\chi_{2,i}\bigl(f(\alpha)\bigr)^2-\chi_{2,i}\bigl(f(\alpha)\bigr)\Bigr). $$ 
\end{lem}

If the characteristic is odd we then use the same arguments as in Section~\ref{sec-repr} to conclude that 
\begin{multline*}  b_{\mu} c_{\nu}|_g = \frac{I}{2^{\abs{\mu}+\abs{\nu}}} \sum_{f \in P_g} \prod_{i=1}^{l(\mu)} \Bigl(\sum_{\alpha \in A(i)}\chi_{2,i}\bigl(f(\alpha)\bigr)+\chi_{2,i}\bigl(f(\alpha)\bigr)^2\Bigr)^{\mu_i} \cdot \\ \cdot\prod_{j=1}^{l(\nu)} \Bigl(\sum_{\alpha \in A(j)}\chi_{2,j}\bigl(f(\alpha)\bigr)-\chi_{2,j}\bigl(f(\alpha)\bigr)^2\Bigr)^{\nu_j}.\end{multline*}
Note that this expression is defined for all $g \geq -1$. It can be decomposed in terms of $u_g$'s (that is, we can find a result corresponding to Lemma~\ref{lem-decomp}) for tuples $(\nn;\rr) \in \mathcal{N}_m$ such that 
\begin{equation} \label{eq-bc-decomp}
\abs{\nn} \leq \abs{\mu}+\abs{\nu}.
\end{equation}

\begin{rmk} The corresponding results clearly hold for elements $(h,f)$ in $P_g$ in even characteristic and the decomposition of $b_{\mu} c_{\nu}|_g$ is independent of characteristic. 
\end{rmk}

\begin{exa} For each $N$ we have the decomposition:
$$b_{[N]}|_g=\frac{1}{2}(u_g^{((N);(2))}+u_g^{((N);(1))}) \; \; \text{and} \; \; c_{[N]}|_g=\frac{1}{2}(u_g^{((N);(2))}-u_g^{((N);(1))}).$$
\end{exa}

\begin{exa} Let us decompose $b_{[1^2]}c_{[2]}|_g$ into $u_g$'s:
\begin{multline*} b_{[1^2]}c_{[2]}|_g=\frac{1}{8}(u_g^{((2,1,1);(2,2,2))}+u_g^{((2,1,1);(2,1,1))}+2u_g^{((2,1);(2,2))}\\-u_g^{((2,1,1);(1,2,2))}-u_g^{((2,1,1);(1,1,1))}-2u_g^{((2,1);(1,2))}).\end{multline*}
In this expression we have removed the $u_g$'s for which $\sum_{i=1}^m r_i n_i$ is odd, since they are always equal to $0$.\end{exa}

\begin{lem} \label{lem-bc-u}
For each $N$, the following information is equivalent:
\begin{itemize}
\item[(1)] all $u_g$'s of degree at most $N$;
\item[(2)] all $b_{\mu} c_{\nu}|_g$ of weight at most $N$.
\end{itemize}
\end{lem}
\begin{pf} From property \eqref{eq-bc-decomp} of the decomposition of $b_{\mu}c_{\nu}|_g$ into $u_g$'s we directly find that if we know $(1)$ we can compute $(2)$.
For the other direction we note on the one hand that 
\begin{equation} \label{eq-create} I \sum_{f \in P_g} \prod_{i=1}^j \bigl(b_i(C_f)-c_i(C_f)\bigr)^{s_i} \bigl(b_i(C_f)+c_i(C_f)\bigr)^{t_i} \end{equation}
can be formulated in terms of $b_{\mu}c_{\nu}|_g$'s of weight at most 
$$S:=\sum_{i=1}^j i \, (s_i+t_i). $$
If we on the other hand decompose \eqref{eq-create} into $u_g$'s we find that there is a unique $u_g$ of degree $S$. The corresponding pair $(\nn;\rr)$ contains, for each $i$, precisely $s_i$ entries of the form $i^1$ and $t_i$ entries of the form $i^2$. Every $u_g$ of degree $S$ can be created in this way and hence if we know $(2)$ we can compute $(1)$. \end{pf}

\begin{rmk} From the definitions of $a_i(C_f)$ and $r_i(C_f)$ we see that knowing $(1)$ and $(2)$ in Lemma \ref{lem-bc-u} is also equivalent to knowing
\begin{itemize}
\item[(3)] all $a_{\lambda}r_{\xi}|_g$ of weight at most $N$,
\end{itemize}
where $a_{\lambda}r_{\xi}|_g$ is defined in the obvious way. Moreover, $a_{\lambda}r_{\xi}|_g=0$ if $\abs{\lambda}$ is odd.
\end{rmk}

\section{Appendix: The stable part of the counts} \label{sec-app2} 
\begin{rmk} All results in this section are independent of characteristic. \end{rmk}

\begin{dfn}[{\cite[Def. 1.2.1, 1.2.2]{DeligneCdWII}}] Let $\mathcal{F}$ be a constructible ($\ell$-adic) sheaf on a scheme $X$ of finite type over $\Z$. 
The sheaf $\mathcal{F}$ is said to be \emph{pure} of weight $m$ if, for every closed point $x$ in $X$ and eigenvalue $\alpha$ of Frobenius $F$ (relative to $k=k(x)$) acting on $\mathcal{F}_{\bar x}$, $\alpha$ is an algebraic integer of weight equal to $m$, i.e., such that all its conjugates have absolute value equal to $q^{m/2}$. The sheaf $\mathcal{F}$ is said to be \emph{mixed} of weight $\leq m$ if there exists a filtration $0=\mathcal{F}_{-1} \subset \mathcal{F}_0 \subset \ldots \subset \mathcal{F}_{m}=\mathcal{F}$ of constructible subsheaves such that, for all $j$, $\mathcal{F}_j/\mathcal{F}_{j-1}$ is pure of weight $j$.
\end{dfn}

\begin{thm}[{\cite[Cor. 3.3.3, 3.3.4]{DeligneCdWII}}] \label{thm-weight} Let $X \xrightarrow{f} \Z$ be a scheme of finite type, and $\mathcal{F}$ a constructible sheaf mixed of weight $\leq m$. Then $R^if_{!}\mathcal{F}$ is mixed of weight $\leq m+i$. Thus, for every finite field $k$, there is a filtration $0=W_{-1} \subset W_0 \subset \ldots \subset W_{i+m}=H^i_{c}(X_{\bar k},\mathcal{F})$ of $\mathrm{Gal}(\bar k /k)$-representations such that, for all $j$, $W_j/W_{j-1}$ is pure of weight $j$.
\end{thm}

\begin{dfn} Let $K_0(\mathsf{Gal}_{k})$ be the Grothendieck group of $\mathrm{Gal}(\bar k /k)$-repre\-sentations. In this category, and with the notation of Theorem \ref{thm-weight}, we have $[H^i_{c}(X_{\bar k},\mathcal{F})]=\sum_{j=0}^{i+m} [W_j/W_{j-1}]$. For any $w \geq 0$, let us define $[H^i_{c}(X_{\bar k},\mathcal{F})]^w:=\sum_{j = w}^{i+m} [W_j/W_{j-1}]$ and $\Eulc^w(X_{\bar k},\mathcal{F}):=\sum_{i \geq 0} (-1)^i [H^i_{c}(X_{\bar k},\mathcal{F})]^w$ in $K_0(\mathsf{Gal}_{k})$. We make the corresponding definition of $\Eulc^w(X_{\overline{\Q}},\mathcal{F})$ in $K_0(\mathsf{Gal}_{\Q})$.  
\end{dfn}

Recall the definition in Section~\ref{sec-coh-et}, for a prime $\ell \nmid q$, of the $\ell$-adic local system $\Vla$ on $\Hh_{g}$. If $\tau$ is the canonical morphism from $\mathcal{H}_g\otimes \bar k$ to $H_g$, we put $\Vla'=\tau_* \Vla$. This is a constructible sheaf pure of weight $\abs{\lambda}$. 

In this section we will see that if $g$ and $w$ are large enough we can compute the trace of Frobenius on $\Eulc^w(\Hh_{g} \otimes \bar k,\Vla)$, which by definition (cf. Section 2 in \cite{jbgvdg}) is equal to $\Eulc^w(H_{g},\Vla')$. We first make the connection to $\s_n$-equivariant counts of points of $H_{g,n}$ explicit.

\begin{lem} \label{lem-loca} Let the symmetric polynomial $\mathbf{s}_{<\lambda>}$ be the Schur polynomial in the symplectic case (see \cite[A.45]{Fulton-Harris}), and $\mathbf{p}_{\lambda}$ the power sum. If $\mathbf{s}_{<\lambda>}= \sum_{\abs{\mu} \leq \abs{\lambda}} m_{\mu} \, \mathbf{p}_{\mu}$ then 
\begin{equation} \label{eq-traceV}
\mathrm{Tr}\bigl(F,\Eulc(\Hh_g \otimes \bar k,\Vla')\bigr)= \sum_{\abs{\mu} \leq \abs{\lambda}} m_{\mu} \, q^{\frac{1}{2}(\abs{\lambda}-\abs{\mu})} \, a_{\mu}|_g.
\end{equation} 
\end{lem}

From Theorems~\ref{thm-rec2} and \ref{thm-rec2even} we see that only the $u_g$'s with all $r_i=2$ have inhomogeneous recurrence relations. Theorem~\ref{thm-chareq} dealt with the homogeneous part of the linear recurrence relations for $a_{\lambda}|_g$. The following lemma, which is a direct consequence of Theorems \ref{thm-rec2}, \ref{thm-rec2even} and \ref{thm-chareq}, deals with the ``inhomogeneities''. 

\begin{lem} \label{lem-part} Denote by $t_{\nn}$ the coefficient of $u_g^{(\nn;(2,\ldots,2))}$ in the decomposition of $a_{\lambda}|_g$ (given in Construction-Lemma~\ref{lem-decomp}). Each value of $\abs{\nn}$ for a pair $(\nn;(2,\ldots,2))$ appearing in this decomposition of $a_{\lambda}|_g$ is at most equal to $\abs{\lambda}/2$. 
Define the polynomial $$f_{\nn}(x):=\bigl(\prod_{i=1}^m (x^{n_i}-1) \bigr)/(x-1).$$ For $g \geq 0$, let $R_{\lambda}(q)|_g$ be the sum, over the pairs $(\nn;(2,\ldots,2))$ that occur in the decomposition of $a_{\lambda}|_g$, of the polynomial quotients of, 
\begin{equation} 
t_{\nn} \, q^{2g+\abs{\nn}} \, J \, (q-1) \, f_{\nn}(q) \quad \; \text{by} \quad \; f_{\nn}(q^2), \end{equation}
which is of degree at most $(\abs{\lambda}+4g-2)/2$. The polynomial $R_{\lambda}(q)|_g$ is a particular solution to the recurrence relation, described in Section~\ref{sec-reca}, for $a_{\lambda}|_g$. 
\end{lem}

Since the power sums form a rational basis of the ring of symmetric polynomials, equation \eqref{eq-traceV} and Theorem \ref{thm-weight} show that $a_{\lambda}|_g$ is of the form $\sum_j z_j \alpha_j$ for a finite set of rational numbers $z_j$ and distinct algebraic integers $\alpha_i$ of weight at most $\abs{\lambda}+4g-2$ (note that $2g-1$ is the dimension of $H_g$). If our base field $k$ is replaced by an extension $k_m$ of degree $m$ then $a_{\lambda}|_g$ is equal to $\sum_j z_j \alpha_j^m$. For $g \geq \abs{\lambda}-1$, the linear recurrence relation for $a_{\lambda}|_g$ (see Section \ref{sec-reca}) shows that it can be written as the particular solution $R_{\lambda}(q)|_g$ plus the homogeneous part, an integer sum of $a_{\lambda}|_{\tilde g}-R_{\lambda}(q)|_{\tilde g}$ for $\tilde g \leq \abs{\lambda}-2$. We then see that if $g \geq \abs{\lambda}-1$ and $w = 5 \, \abs{\lambda}-9$, the homogeneous part of the solution to the linear recurrence relation for $a_{\lambda}|_g$ does not contribute to $\mathrm{Tr}\bigl(F,\Eulc^w(\Hh_g \otimes \bar k,\Vla')\bigr)$. To conclude this we used the fact that $\sum_i z_i \alpha_i^m=0$ for all $m$ implies that $z_i=0$ for all $i$, where the $z_i$ and $\alpha_i$ are complex numbers and the $\alpha_i$ are distinct and nonzero. We can now summarize using Theorem 3.2 in \cite{Jonas2}. 

\begin{dfn}
For a polynomial $f(x)=\sum_i f_i x^i$ put $f^w(x):=\sum_{i \geq w} f_i x^i$. 
\end{dfn}

\begin{thm} \label{thm-stable} 
Let $\mathbf{q}$ denote the class of $\Ql(-1)$ in $K_0(\mathsf{Gal}_{\Q})$.
For  $g \geq \abs{\lambda}-1$ and $w = 5 \, \abs{\lambda}-9$ we have an equality in $K_0(\mathsf{Gal}_{\Q})$,
$$\Eulc^w(\Hh_{g} \otimes \overline{\Q},\Vla)= \sum_{\abs{\mu} \leq \abs{\lambda}} m_{\mu} \, \mathbf{q}^{\frac{1}{2}(\abs{\lambda}-\abs{\mu})} \, R^{w-\abs{\lambda}+\abs{\mu}}_{\mu}\bigl(\mathbf{q}\bigr)|_g.$$
\end{thm}

\begin{exa} In the case $\lambda=(4,2,2)$, for $w=31$ and $g \geq 7$, we find that $\mathrm{Tr}\bigl(F,\Eulc^w(\Hh_{g} \otimes \bar k,\Vla)\bigr)$ is equal to $f_g^w(q)$, where $f_g$ is the polynomial quotient of $q^{2g+4}(3q^2+3q+2)$ by $(q^2+1)^2(q+1)^3$. 
\end{exa}

\begin{rmk} By Poincar\'e duality (cf. Section 2 in \cite{jbgvdg}) we find that there is a filtration $0=W'_{i+\abs{\lambda}-1} \subset W'_{i+\abs{\lambda}} \subset \ldots \subset W'_{2(2g-1+\abs{\lambda})}=H^i(\Hh_g \otimes \bar k,\Vla)$ of $\mathrm{Gal}(\bar k /k)$-representations such that $W'_j/W'_{j-1}$ is pure of weight $j$. Let us define $[H^i(\Hh_g \otimes \bar k,\Vla)]_w:=\sum_{j=i+\abs{\lambda}}^{w} [W'_j/W'_{j-1}]$ and $\Eul_{w}(\Hh_g \otimes \bar k,\Vla):=\sum_{i \geq 0} (-1)^i [H^i(\Hh_g \otimes \bar k,\Vla)]_w$ in $K_0(\mathsf{Gal}_{k})$ and similarily $\Eul_{w}(\Hh_{g} \otimes \overline{\Q},\Vla)$. Theorem \ref{thm-stable} shows that, for $g \geq \tilde g \geq \abs{\lambda}-1$ and $w = 4 \tilde g -3\,\abs{\lambda}+7$, one has that $\Eul_{w}(\Hh_{g} \otimes \overline{\Q},\Vla)$ is \emph{stable}, in the sense that it is independent of $g$.
\end{rmk}

Computations for $\lambda$'s of low weight lead us to make a conjecture, which is true for $\abs{\lambda} \leq 30$. 
\begin{conj} For  $g \geq \abs{\lambda}-1$ and $w = 5 \, \abs{\lambda}-9$, we have $\Eul_{c}^w(\Hh_{g} \otimes \overline{\Q},\Vla)=0$ for all $\lambda$ such that $\lambda_1 > \abs{\lambda}/2$. 
\end{conj}

\bibliographystyle{plain}

\end{document}